\documentclass[a4paper,10pt]{article}

\usepackage[utf8]{inputenc}
\usepackage[english]{babel}
\usepackage[babel]{csquotes}
\usepackage{amsthm, amsmath, amssymb,amsfonts,latexsym,soul}
\usepackage{mathtools} 
\usepackage[inline]{enumitem}  
\usepackage{fullpage}  
\usepackage{bbm}       
\usepackage{mathrsfs}  
\usepackage{mhchem}  

    \usepackage[numbers,sort]{natbib}

\usepackage[colorlinks]{hyperref}

\newcommand{\Sp}{\mathcal{X}}
\newcommand{\Spacr}{\Sp_{\mathrm{ACR}}}
\newcommand{\piacr}{\pi_{\mathrm{ACR}}}
\newcommand{\C}{\mathcal{C}}
\newcommand{\R}{\mathcal{R}}
\newcommand{\G}{\mathcal{G}}
\newcommand{\K}{\mathcal{K}}
\newcommand{\Sy}{\mathscr{S}}
\newcommand{\NN}{\mathbb{N}}
\newcommand{\RR}{\mathbb{R}}
\newcommand{\Pois}{\mathrm{Pois}}
\newcommand{\pr}[1]{\left(#1\right)}
\newcommand{\abs}[1]{\left|#1\right|}
\newcommand{\prq}[1]{\left[#1\right]}

\newcommand{\Xd}{X_{\mathrm{disc}}}
\newcommand{\Xacr}{X_{\mathrm{ACR}}}
\newcommand{\oXd}{\overline{X}_{\mathrm{disc}}}
\newcommand{\Xc}{X_{\mathrm{cont}}}
\newcommand{\oXc}{\overline{X}_{\mathrm{cont}}}
\newcommand{\fastR}{\widetilde{\R}}
\newcommand{\fastK}{\widetilde{\K}}
\newcommand{\Syl}{\widetilde{\Sy}}
\newcommand{\floor}[1]{\lfloor #1 \rfloor}

\DeclareMathOperator{\spann}{span}

\binoppenalty=\maxdimen
\relpenalty=\maxdimen

\newtheorem{theorem}{Theorem}[section]
\newtheorem{corollary}[theorem]{Corollary}

\theoremstyle{remark}
\newtheorem{remark}{Remark}[section]
\newtheorem{example}{Example}[section]

\theoremstyle{definition}
\newtheorem{definition}{Definition}[section]
\newtheorem{assumption}{Assumption}

\begin{document}

\title{Finite time distributions of stochastically modeled chemical systems with absolute concentration robustness}

\author{David F.\ Anderson\footnotemark[1] \and Daniele Cappelletti\footnotemark[2] \and Thomas G.\ Kurtz\footnotemark[3]}

\footnotetext[1]{Department of Mathematics, University of Wisconsin-Madison, email: \texttt{anderson@math.wisc.edu}. Supported by NSF-DMS-1318832 and Army Research Office grant W911NF-14-1-0401.}
\footnotetext[2]{\emph{Corresponding author}. Department of Mathematical Sciences, University of Copenhagen, email: \texttt{cappelletti@math.wisc.edu}. Supported by the Danish Research Councils.}
\footnotetext[3]{Department of Mathematics, University of Wisconsin-Madison, email: \texttt{kurtz@math.wisc.edu}. Supported by NSF-DMS-11-06424.}
\maketitle
 
\begin{abstract} 
 Recent research in both the experimental and mathematical communities has focused on biochemical interaction systems that satisfy an ``absolute concentration robustness'' (ACR) property.   The ACR property was first discovered experimentally when, in a number of different systems, the concentrations of key system components at equilibrium were observed to be robust to the total concentration levels of the system.  Followup mathematical work focused on deterministic models of biochemical systems and demonstrated how chemical reaction network theory can be utilized to explain this robustness.   Later mathematical work  focused on the behavior of this same class of reaction networks, though under the assumption that the dynamics were stochastic.  Under the stochastic assumption, it was proven that the system will undergo an extinction event with a probability of one so long as the system is conservative, showing starkly different long-time behavior than in the deterministic setting.  Here we consider a general class of stochastic  models that intersects with the class of ACR systems studied previously.  We consider a specific system scaling over compact time intervals and prove that in a  limit of  this scaling the distribution of the  abundances of the ACR species converges to a certain product-form Poisson distribution whose mean is the ACR value of the deterministic model. This result is in agreement with recent conjectures pertaining to the behavior of  ACR networks endowed with stochastic kinetics, and helps to resolve the conflicting theoretical results pertaining to deterministic and stochastic models in this setting.
\end{abstract}

\section{Introduction}

 Biochemical reaction networks are often quite complex and computationally intractable. It is therefore important to develop mathematical techniques that relate simple graphical features of the reaction network, which are easy to check, to the qualitative dynamics of the underlying mathematical model. This approach dates back to at least \cite{feinberg:lectures,feinberg_horn,horn:general_mass_action}, where certain graphical characteristics of networks were shown to ensure uniqueness and local asymptotic stability of the steady states for deterministically modeled complex-balanced systems.  
 
 In this context of relating graphical and dynamical features of  models, Shinar and Feinberg  provided graphical conditions  that imply certain species satisfy an  \emph{absolute concentration robustness} (ACR) property for the associated deterministically modeled system \cite{shinar:structural}.
 A species is said to possess ACR if for a fixed choice of system parameters its concentration is the same at any positive equilibrium point of the deterministically modeled system. Such a feature has been observed experimentally in several important biochemical reaction networks, including signal transduction cascades and gene regulatory networks \cite{blanchini:structurally,shinar:structural,shinar:sensitivity,alon:robustness,dexter:dimerization,gunawardena:models}. The ACR property provides useful information on the system dynamics since it indicates a predictable fixed response regardless of  changes in the environment. Followup research pertaining to deterministically modeled systems with ACR species can be found in \cite{Gun2014}.
  
 Stochastically modeled systems satisfying essentially the same graphical conditions as those detailed  in \cite{shinar:structural} were considered by Anderson, Enciso, and Johnston in \cite{anderson:ACR}.  There it was shown that ACR systems in a particular family, if stochastically modeled, undergo an extinction event with a probability of one, so long as the system is conservative (i.e.\ there is a positive linear combination of the species that is invariant to occurrences of reactions).
 Such a result can be considered  an example of a discrepancy between the limiting behaviour of a deterministic system and the limiting behaviour of the corresponding stochastic system, with one modeling choice predicting a form of long-term stability and the other predicting long-term instability. However, in \cite{anderson:ACR} it is pointed out that the extinction event is typically a rare event on reasonable timeframes and that useful information pertaining to the behaviour of stochastically modeled ACR systems could be had by better understanding the dynamics of the system on compact time intervals or via the quasi-stationary distribution.  It is conjectured in \cite{anderson:ACR} that the distribution of the ACR species will be approximately Poisson in either case.  Both a  simple example pertaining to a model of protein interactions and a numerical analysis of the two-component EnvZ/OmpR signaling system in \textit{Escherichia coli} provide evidence in favor of  the conjecture \cite{anderson:ACR}. 
 
 In this paper we provide an asymptotic result for the stochastic models of a class of reaction systems that overlaps with ACR reaction systems. In particular, we consider a multiscale setting in which the abundances of a subset of the ACR species are of order $O(1)$, while the abundances of other species are of order $O(N)$.  We then scale the rate constants in a particular way and let $N$ go to infinity. Under this limit, we prove that on compact time intervals the ACR species whose abundance is of order $O(1)$ behave in the way conjectured in \cite{anderson:ACR}. Namely, the distribution of their abundances is well approximated by a product-form Poisson distribution whose parameter is given by the ACR equilibrium value of the associated deterministically modeled system. Thus, the results presented here link the qualitative behaviours of the deterministic and stochastic models. Furthermore, the result fully explains the outcome of the numerical analysis of the EnvZ/OmpR signaling system performed in \cite{anderson:ACR}.  
 
One key observation we utilize in our proofs, and the basis of our Assumption \ref{ass:complex_balanced},   is that a certain sub-reaction network consisting of only the $O(1)$ species
 is often weakly reversible and deficiency zero.  This fact, together with the results of \cite{anderson:product-form}, allows us to characterize the marginal distribution of the ACR species as approximately Poisson.  We then show that the approximation becomes precise in the limit as $N\to \infty$.
 
 We end this section with two instructive examples that demonstrate our main results.
 
 \begin{example}\label{ex:simple}
  Consider the deterministically modeled system with reaction network
  \begin{equation}\label{eq:ex_simple}
    A + B \ce{->[\kappa_1]} 2B,    \qquad B \ce{->[\kappa_2]} A
  \end{equation}
  and  mass action kinetics (see  \eqref{eq:dma1} and \eqref{eq:dma2}). 
 The species $A$ exhibits ACR since at each positive equilibrium the concentration of $A$ is $\kappa_2/\kappa_1$, regardless of the concentration of the species $B$ \cite{anderson:ACR,shinar:structural}. 
 
Now consider a sequence $\{X^N\}_{N\in\NN}$ of continuous time Markov chain models for \eqref{eq:ex_simple}, in which the counts of species $A$ and $B$ at time $t$ are given by $X^N_1(t)$ and $X^N_2(t)$, respectively. We suppose the initial conditions are such that $X_1^N(0)$ is a bounded sequence and $N^{-1}X^N_2(0)$ converges to a positive real number, as $N\to \infty$. This choice of initial conditions corresponds to an experiment where the abundance of the molecules of $B$ is increased, while the magnitude of the count of $A$ is maintened. Our goal will be to understand the limiting behaviour of $X^N_1$ for $N$ going to infinity.  

For the sake of intuition, note that  a sub-reaction network for $A$ is 
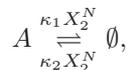
\[
	A\overset{\kappa_1 X_2^N}{\underset{\kappa_2X_2^N}{\rightleftharpoons}} \emptyset,
\]
which views $X_2^N$ as simply modulating the speed of the two reactions.  Since it is well known that in the case of the above reaction network with $X_2^N$ fixed, $X_1^N$ will have a stationary distribution that is Poisson with parameter $\kappa_2/\kappa_1$ \cite{anderson:product-form}, it is intuitively clear that for finite $t$, $X_1^N(t)$ should be approximately Poisson if $X_2^N$ is large and fluctuates little.

The above argument will be made precise.  Let $J$ be a random variable having a Poisson distribution with parameter $\kappa_2/\kappa_1$. Corollary \ref{cor:ACR_equilibrium_point} will allow us to conclude that for any function $g:\NN\to\RR$ with at most polynomial growth rate, the distribution of $g(X^N_1(t))$ converges on average to the distribution of $g(J)$ for $N$ going to infinity. Specifically, for any real positive $T$, we have 
 \[
 	\sup_{t\in[0,T]}\abs{\int_0^t\Big(g(X_1^N(s))-E[g(J)]\Big)ds}\xrightarrow[N\to\infty]{}0
 \]
 in probability.  Further, if $g$ is bounded, then for any $\delta>0$
 \begin{equation}\label{eq:56865}
 	\sup_{t\in[\delta,T]}\abs{E\Big[g(X_1^N(t))\Big]-E[g(J)]}\xrightarrow[N\to\infty]{}0.
 \end{equation}
  An immediate consequence of \eqref{eq:56865} is that $X_1^N(t)$ converges in distribution to $J$, which has a Poisson distribution with mean given by the ACR value $\kappa_2/\kappa_1$.  This consequence follows from choosing as $g$ an indicator function. 
See Example \ref{ex:ACR_equilibrium_point} for more details.
\hfill $\square$
 \end{example}
 \begin{example}\label{ex:non_mass action}
  Our results are also applicable to models that  do not utilize mass action kinetics.  Consider the stochastic reaction system
  
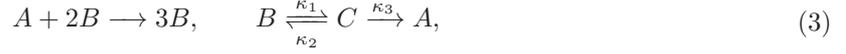
\begin{gather}\label{eq:ex_not_mass_action}
   \begin{split}
    A + 2B \xrightarrow{\phantom{\kappa_3}} 3B, \qquad 
    B\ce{<=>[\kappa_1][\kappa_2]} C \xrightarrow{\kappa_3} A,
   \end{split}
  \end{gather}
  where the  rate of the reaction $A+2B\to 3B$ is given by
  \begin{equation}
  \label{eq:0987}
  \lambda(x)=\kappa_0\frac{x_1x_2(x_2-1)}{1+x_2},
 \end{equation}
  and where $\kappa_0\in\RR_{>0}$. The rate \eqref{eq:0987} corresponds to an inhibitory effect of the molecules of $B$ on the production of $B$ itself. If we consider a sequence of such models in which the counts of $B$ and $C$ go to infinity, then the limiting behaviour of the model \eqref{eq:ex_not_mass_action}-\eqref{eq:0987}  coincides with the limiting behaviour of the process associated with the reaction network
  
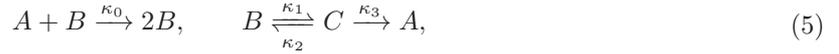
\begin{gather}
  \label{eq:68750}
  \begin{split}
   A + B \xrightarrow{\kappa_0} 2B, \qquad 
   B\ce{<=>[\kappa_1][\kappa_2]} C \xrightarrow{\kappa_3} A,
   \end{split}
  \end{gather}
  endowed with mass-action kinetics.  Due to \cite{shinar:structural}, the  reaction system \eqref{eq:68750} exhibits ACR in the species $A$, when deterministically modeled. Let $q\in\RR_{>0}$ be the ACR value for species $A$.  Due to the connection between the models \eqref{eq:ex_not_mass_action}-\eqref{eq:0987} and \eqref{eq:68750}, we anticipate the value $q$ will play a role in the limiting behaviour of species $A$ of \eqref{eq:ex_not_mass_action}-\eqref{eq:0987}.
  
We therefore denote by $\{X^N\}_{N\in\NN}$ a sequence of stochastic processes modeled according to  \eqref{eq:ex_not_mass_action}-\eqref{eq:0987},  with $X^N_1(t)$, $X^N_2(t)$, and $X^N_3(t)$ being the counts at time $t$ for the species  $A$, $B$, and $C$, respectively. We suppose that $X^N_1(0)$ is a bounded sequence and $N^{-1}(X^N_2(0),X^N_3(0))$ converges to a point in $\RR^2_{>0}$, as $N\to \infty$.
Our aim is to understand the limiting behaviour of  $X^N_1$ as  $N$ goes to infinity.
  
  Let $J$ be a Poisson distribution with parameter $q$, the ACR value for species $A$ when \eqref{eq:68750} is modeled deterministically. By Corollary \ref{cor:ACR_basin}, there exits a function $g^*:\RR_{>0}\to\RR_{>0}$ with $g^*(s)\to E[g(J)]$, as $s\to\infty$, such that for any positive real $T$
  $$\sup_{t\in[0,T]}\abs{\int_0^{t}\Big(g(X_1^N(s))-g^*(s)\Big)ds}\xrightarrow[N\to\infty]{}0$$
  in probability.  Moreover, and as in Example \ref{ex:simple}, we have that $X_1^N(t)$ converges in distribution to a Poisson random variable whose mean is related to the ACR value.
  See Example \ref{ex:ACR_basin} for more details.
   \hfill $\square$  
 \end{example}

\section{Necessary Background and Notation}
\label{sec:background}

 We denote the natural numbers including 0 by $\NN$, that is $\mathbb N = \{0,1,2,\dots\}$. For any real vector $v$, we denote its $i$th entry by $v_i$. We will write $v>0$ if every entry of $v$ is strictly positive. We denote by $[v]$ the vector of the floor functions of the entries of $v$; that is, $[v]_i = \floor{v_i}$. For any real vector $\alpha$ of the same size as $v$, and for $N>0$, we denote by $N^{\alpha}v$ the  vector satisfying
 \[
 	\pr{N^{\alpha}v}_i=N^{\alpha_i}v_i.
 \]
 We will denote by $\|v\|$ the euclidean norm of the vector, by $\|v\|_1$ its \emph{$L^1$-norm} and by $\|v\|_\infty$ its \emph{$L^\infty$-norm}, that is
 \[
 	\|v\| = \sqrt{\sum_i v_i^2}, \qquad 
	 \|v\|_1=\sum_i|v_i|, \quad\text{and}\quad \|v\|_\infty=\max_i|v_i|.
\]
 For two vectors $v$ and $w$ of the same dimension, we write $v<w$, $v\leq w$, $v>w$ or $v\geq w$ if the inequality holds component-wise. Furthermore, for any set $A$ we will indicate by $|A|$ its cardinality and by $\mathbbm{1}_A$ its indicator function.
 Finally, for $a,b \in \mathbb{R}$, we denote $a\wedge b = \min\{a,b\}$ and $a\vee b = \max\{a,b\}$.
 
 We say that a function $g\colon \RR^n\to\RR$ has at most \emph{polynomial growth rate} if there exists a multivariate polynomial $p\colon \RR^n\to\RR$ such that
 $$\limsup_{\|x\|\to\infty}\frac{g(x)}{p(x)}=0$$

 Here we give some basic definitions from chemical reaction network theory, see for example \cite{erdi:mathematical_models,feinberg:lectures} for a more detailed introduction.
 
 A \emph{reaction network} is a triple $\G=(\Sp,\C,\R)$. $\Sp$ is a finite non-empty ordered set of symbols, referred to as \emph{species}, and $\C$ is a finite non-empty ordered set of linear combinations of species with non-negative integer coefficients, referred to as \emph{complexes}. Any species $S_i\in\Sp$ can be identified with the vector $e_i\in\RR^{|\Sp|}$, whose $i$th entry is 1 and whose other entries are zero. Therefore, any complex $y\in\C$ will be identified with a vector in $\RR^{|\Sp|}$ that is linear combination of the vectors $e_i$. Finally, $\R$ is a non-empty ordered subset of $\C\times\C$, whose elements are called \emph{reactions}, such that for any $y\in\C$, $(y,y)\notin\R$. Following the common notation, we will denote any element $(y_r,y_r')\in\R$ by $y_r\to y_r'\in\R$, in which case we then call $y_r$ the \emph{source complex} and $y_r'$ the \emph{product complex} of that reaction. It is possible that a complex $y\in\C$ is the source (product) complex of different reactions, and that it is both the source complex of one reaction and the product complex of another reaction. It is commonly required that every species $S\in\Sp$ appears in at least one complex, and that every complex $y\in\C$ appears as an element in at least one reaction. It is possible to associate a directed graph to $\G$, where the set of nodes is the set of complexes $\C$ and the arrows are given by the reactions $y_r\to y_r'\in\R$. If the graph is such that for any directed path from $y$ to $y'$ there exists a directed path from $y'$ to $y$, then $\G$ is \emph{weakly reversible}.
 For the $r$th reaction, $y_r\to y_r'$, we denote by $\xi_r=y_r' - y_r$ the corresponding \emph{reaction vector}. We write $S_i \in\xi_r$ ($S_i \in y_r$) if $\xi_{ri} \ne 0$ ($y_{ri} \ne 0$). For any species $S\in\Sp$, let 
 \begin{equation}\label{eq:78965678}
 	\R_S=\{y_r\to y_r'\in\R\,:\, S\in\xi_r\},
\end{equation}
 the set of reactions that change the amount of species $S$.
 
To each reaction $y_r\to y_r'\in\R$, we can associate a function $\lambda_r:\RR^{|\Sp|}_{\geq0}\to\RR_{\geq0}$. The set consisting of  these function $\K=\{\lambda_r\}_{y_r\to y_r'\in\R}$ is referred to as the \emph{kinetics}, and the functions $\lambda_r$ are called \emph{rate functions}, or \emph{intensity functions}, or \emph{propensity functions}. The pair $\Sy=(\G,\K)$ is a \emph{reaction system}, which can be stochastically or deterministically modeled, as explained below.
 
 In a \emph{stochastically modeled reaction system} $\Sy=(\G,\K)$, the counts  of molecules of the different chemical species are considered, and the  counts at time $t$ form a vector $X(t)\in\NN^{|\Sp|}$. The evolution in time of the vector $X(t)$ follows a continuous time Markov chain, where in each state $x\in\NN^{|\Sp|}$ the obtainable states are $\{x+\xi_r\}$ and transition rates are given by $\{\lambda_r(x)\}$, with $y_r\to y_r'$ varying in $\R$. If at time $t^*$ the reaction $y_r\to y_r'$ occurs, then we have
 $$X(t^*)=X(t^*-)+\xi_r,$$
 where $X(t^*-)$ denotes the previous  state. To stick with the physical meaning of the reactions, we require that the kinetics is such that for any reaction $y_r\to y_r'\in\R$ we have $\lambda_{y_r\to y_r'}(x)>0$ only if $x\geq y_r$. This condition prevents the number of molecules present from becoming negative.  Moreover, in this setting, we are only interested in the values of $\lambda_{y_r\to y_r'}(x)$, when $x\in\NN^{|\Sp|}$, therefore the domain of the rate function can be restricted to $\NN^{|\Sp|}$. Following the terminology utilized in \cite{anderson:design,AK2015,ethier_kurtz,kurtz:strong}, we can write
 $$X(t)=X(0)+\sum_{y_r\to y_r'\in\R} Y_r\pr{\int_0^t\lambda_r(X(s))ds} \xi_r,$$
 where the $Y_r$ are i.i.d.~unit rate Poisson processes.
 For any two states $x,z\in\NN^{|\Sp|}$, we say that a state $z$ is \emph{obtainable} from $x$ if there exists a sequence of reactions $(y_{r_i}\to y_{r_i}')_{i=1}^m$ such that
 $$z=x+\sum_{i=1}^m\xi_{r_i}, \quad \text{ with }\quad  \lambda_{y_{r_i}\to y_{r_i}'}\pr{x+\sum_{j=1}^{i-1}\xi_{r_j}}>0,$$
 for all $i \in \{1,\dots,m\}$.
 We further say that $(\G,\K)$ is \emph{irreducible} if for any two states $x,z\in\NN^{|\Sp|}$, $z$ is obtainable from $x$ and $x$ is obtainable from $z$. See \cite{craciun:dynamical} for more on irreducible reaction networks and for sufficient conditions implying irreducibility.
 A popular choice of kinetics for stochastic reaction systems is given by \emph{stochastic mass action kinetics}, defined by
 $$\lambda_r(x)=\kappa_r\frac{x!}{(x-y_r)!}\mathbbm{1}_{\{x\geq y_r\}},$$
 where $\kappa_r\in\RR_{>0}$ are called \emph{rate constants} and for any vector $v\in\NN^{m}$, $v!$ is defined by
 $v!=\prod_{i=1}^m v_i!,$
 with the convention $0!=1$. This kinetics is related to the assumption that the system is well-stirred, so the propensity of each reaction is proportional to the number of possible sets of molecules that can give rise to an occurrence of the reaction. A stochastic reaction system endowed with stochastic mass action kinetics is referred to as \emph{stochastic mass action system}, and will be denoted $\Sy=(\G,\kappa)$. Note that the property of irreducibility of a mass action system does not depend on the particular choice of rate constants: indeed, in mass action systems a rate $\lambda_r(x)$ is strictly positive if and only if $x\geq y_r$. 
 
 In a \emph{deterministically modeled reaction system} $\Sy=(\G,\K)$, the concentrations of the different chemical species are considered, and the  concentrations at time $t$ form a vector  $z(t)\in\RR_{\geq0}^{|\Sp|}$.   The evolution in time of the vector $z(t)$ obeys the ordinary differential equation (ODE)
  \begin{equation}
  \label{eq:dma1}
  	z'(t) = \sum_{y_r\to y_r'\in\R}\xi_r \lambda_r(z(t)).
  \end{equation}
 As in the stochastic case, we put a restriction on the kinetics and require that for any $y_r\to y_r'\in\R$ we have $\lambda_{y_r\to y_r'}(x)>0$ only if $x_i>0$ whenever $S_i\in y_r$. This condition means that a reaction cannot take place if some necessary chemical species is missing, and it guarantees that the vector $z(t)$ will remain non-negative.  \emph{Deterministic mass action kinetics} is given by
 \begin{equation}
 \label{eq:dma2}
 	\lambda_r(x)=\kappa_rx^{y_r},
\end{equation}
 where $\kappa_r\in\RR_{>0}$ are called \emph{rate constants} and for any two vectors $v,w\in\NN^{m}$, $v^w$ is defined by
 $v^w=\prod_{i=1}^m v_i^{w_i},$
 with the convention $0^0=1$. Thus,   the rate of each reaction is proportional to the products of the concentrations of the species appearing in the source complex, according to multiplicity.  As in the stochastic case, this kinetics is chosen for well-stirred systems.   A deterministic reaction system with deterministic mass action kinetics is termed  a \emph{deterministic mass action system}, and will be denoted by $\Sy = (\G,\kappa)$.
 
  A fruitful notion in chemical reaction network theory, and one that will play a role in the present work, is that of a \emph{complex balanced equilibrium}, which is a positive equilibrium point $c$ of a deterministic mass action system satisfying
 $$\sum_{\substack{y_r\to y_r'\in\R\\y_r=y\phantom{'}}}\kappa_rc^{y_r}=\sum_{\substack{y_r\to y_r'\in\R\\y_r'=y}}\kappa_rc^{y_r}\quad\text{for each }y\in\C,$$
where the sum on the left, respectively right, is over those reactions for which $y$ is the source, respectively product, complex.  We say that a deterministic mass action system is \emph{complex balanced} if there exists at least one positive equilibrium point, and if every positive equilibrium point is a complex balanced equilibrium. 
 
 We extend the definition of complex balanced to the stochastic setting by saying that a stochastic mass action system $(\G,\kappa)$ is complex balanced if the deterministic mass action system $(\G,\kappa)$ is complex balanced. We may therefore refer to complex balanced mass action systems without specifying whether they are stochastically or deterministically modeled. In the same fashion, whenever we refer to an equilibrium point of a reaction system, we implicitly assume it is an equilibrium point for the deterministically modeled system. It is interesting to point out that complex balanced stochastic mass action systems can be fully characterised by properties of their stationary distributions \cite{cappelletti:complex_balanced}.
 
 It is worth noting that under the assumptions detailed above for both  deterministic and stochastic reaction systems, the evolution of the amounts of species present is restricted to
 $$X(t)\in \pr{X(0)+\spann\{\xi_r\}_{y_r\to y_r'\in\R}}\cap \NN^{|\Sp|}\quad\text{and}\quad z(t)\in \pr{z(0)+\spann\{\xi_r\}_{y_r\to y_r'\in\R}}\cap \RR_{\geq0}^{|\Sp|},$$
 regardless of the choice of kinetics $\K$. The sets $\pr{v+\spann\{\xi_r\}_{y_r\to y_r'\in\R}}$ with $v \in \RR^{|\Sp|}$ are called the \emph{stoichiometric compatibility classes} of $\G$, and the sets
 \[
 	 \pr{v+\spann\{\xi_r\}_{y_r\to y_r'\in\R}}\cap \RR_{\geq0}^{|\Sp|} \quad \text{ and }\quad \pr{v+\spann\{\xi_r\}_{y_r\to y_r'\in\R}}\cap \RR_{>0}^{|\Sp|}
	 \]
	    are called the \emph{non-negative stoichiometric compatibility classes} and \emph{positive stoichiometric compatibility classes} of $\G$. Any vector $T\in\RR^{|\Sp|}$ that is orthogonal to the stoichiometric compatibility classes of $\G$ is a \emph{conservation law} for $\G$, and if there exists a positive conservation law for $\G$, then $\G$ is called \emph{conservative}.
 
 Let 
 $s=\dim\pr{\spann\{\xi_r\}_{y_r\to y_r'\in\R}}.$
 We define the \emph{deficiency} of $\G$ as 
 $\delta=|\C|-\ell-s,$
 where $\ell$ is the number of connected components of the directed graph associated with $\G$. We end this section by stating some classical results that can be found in  \cite{feinberg:lectures,feinberg_horn,horn:general_mass_action}, which connect graphical and dynamical features of the deterministic mass action systems and will be of use to us.
  \begin{theorem}\label{thm:complex_balanced}
  If a deterministic mass action system $\Sy=(\G,\kappa)$ possesses a complex balanced equilibrium, then $\Sy$ is complex balanced and $\G$ is weakly reversible. Moreover, there exists exactly one complex balanced equilibrium in every positive stoichiometric compatibility class, and it is locally asymptotically stable relative to its positive stoichiometric compatibility class. 
 \end{theorem}
 \begin{theorem}\label{thm:deficiency_zero_iff}
  If $\G$ is weakly reversible and has deficiency 0, then for any choice of rate constants the deterministic mass action  system $\Sy=(\G,\kappa)$ is complex balanced.
 \end{theorem}
   
\section{The multiscale setting, assumptions and main results}\label{sec:general_setting}

We begin by motivating the scaling presented below.  Our goal is to study the behavior of the distribution of  ACR species in the limit as total abundances go to infinity.  
As the equilibrium value of the ACR species is independent of total abundances, we will assume a partition in the set of species: some of them will be allowed to take arbitrarily  large abundances in their initial conditions, while  the initial conditions for the others (i.e.~the ACR species) will be bounded.

 Formally, denote by $\K^N$ a sequence of stochastic kinetics for $\G$, with $N\in\NN_{>0}$, and let $X^N(t)$ be the sequence of stochastic processes associated with the system $(\G,\K^N)$. Assume that there exists a vector $\alpha\in\{0,1\}^{|\Sp|}$ such that
 \begin{align}
 \label{eq:0987897}
 \lim_{N\to\infty}N^{-\alpha}X^N(0)=X_0>0.
 \end{align}
 The condition \eqref{eq:0987897} implies  a partition of the set of species $\Sp$ in two sets, the \emph{discrete species} (denoted by $\Sp_d$) and the \emph{continuous species} (denoted by $\Sp_c$),  
 \begin{itemize}
 \item $S_i \in \Sp_d$ if  $\alpha_i = 0$, in which case $X_i^N(0)=O(1)$;
 \item $S_i\in\Sp_c$ if  $\alpha_i = 1$, in which case $X_i^N(0)=O(N)$.
 \end{itemize}
 Let
 $$\pi_d\colon \RR^{|\Sp|}\to\RR^{|\Sp_d|}\quad\text{and}\quad\pi_c\colon \RR^{|\Sp|}\to\RR^{|\Sp_c|}$$
 be the projections onto the discrete and continuous species, respectively, and define
 $$\Xd^N(t)=\pi_d(X^N(t))\quad\text{and}\quad\Xc^N(t)=\pi_c(X^N(t)).$$
 For convenience, we will sometimes consider the rate functions as functions from $\NN^{|\Sp_d|}\times\NN^{|\Sp_c|}$, and write $\lambda_r(v,w)$, where $v$ and $w$ denote the amounts of the discrete and continuous species, respectively. For any reaction $y_r\to y_r'\in\R$, define
 $$\beta_r=\max_{S_i\in y_r}{\alpha_i}.$$
 Note that $\beta_r\in\{0,1\}$. We assume that 
 \begin{equation}\label{eq:lambda_convergence}
  \lim_{N\to\infty}N^{-\beta_r}\lambda^N_r(v,[Nw])=\lambda_r(v,w)
 \end{equation}
 uniformly on the compact sets of $\NN^{|\Sp_d|}\times\RR_{\geq0}^{|\Sp_c|}$, where the functions $\lambda_r$ are non-zero and locally Lipschitz, with domain $\NN^{|\Sp_d|}\times\RR_{\geq0}^{|\Sp_c|}$. We denote by $\K$ the kinetics given by the limiting functions $\lambda_r$. 
 The above setting is a particular case of the one studied in \cite{ball:asymptotic,kurtz:rescale,popovic:rescale}. 
 \begin{remark}\label{rem:rate_constants_convergence}
  If the kinetics $\K^N$ are stochastic mass action kinetics for all $N\in\NN_{>0}$, then by \eqref{eq:lambda_convergence} the sequence
  $$N^{-\beta_r}\lambda^N_r(v,[Nw])=N^{-\beta_r}\kappa^N_r\frac{v!}{(v-\pi_d(y_r))!}\frac{[Nw]!}{([Nw]-\pi_c(y_r))!}$$
  converges to a positive number for some $(v,w)\in\NN^{|\Sp_d|}\times\RR_{\geq0}^{|\Sp_c|}$. Since the sequence
  $$N^{-\|\pi_c(y_r)\|_1}\frac{[Nw]!}{([Nw]-\pi_c(y_r))!}$$
  converges to $w^{y_r}$, for any $r$ we have
  \begin{equation}\label{eq:rate_constants_convergence}
   \lim_{N\to\infty}N^{-\beta_r+\|\pi_c(y_r)\|_1}\kappa^N_r=\kappa_r,
  \end{equation}
  for some positive constant $\kappa_r$.
 \end{remark}

 \subsection{Assumptions}
 
While the assumptions we detail in this section are technical in nature, they essentially ensure three very natural conditions. Assumption \ref{ass:discrete=fast} ensures that each discrete species is  produced and consumed at high rate. Assumptions \ref{ass:existence_limit} and \ref{ass:moments_bound} ensure that the processes do not explode in finite time.  Finally, Assumption \ref{ass:complex_balanced} requires that the reduced system obtained by the deletion of the high abundance species is complex-balanced.
 
 \begin{assumption}\label{ass:discrete=fast}
 For any $S\in\Sp_d$, there exists at least one reaction $y_r\to y_r'\in\R_S$ such that $\beta_r=1$ (i.e.\ the species $S$ is \emph{fast} consumed or produced, recall $\R_S$ defined in \eqref{eq:78965678}).
 \end{assumption}

 In order to motivate and explain the above assumption  we consider  sequences of processes satisfying the network structures of Examples \ref{ex:simple} and \ref{ex:non_mass action}. 

Consider first the network of Example \ref{ex:simple} with rate constants $\kappa_1$ and $\kappa_2$, and suppose that the total initial abundance of the system (i.e.~the sum of the abundances of species $A$ and $B$) is large, and that the system is near the known ACR equilibrium, in which case $X_1(0) \approx q = \kappa_2/\kappa_1$.  Specifically, we suppose that $X^N_1(0)+X^N_2(0)=N$ for some large $N \in \NN$ and that $X^N_1(0) = O(1)$ in $N$, in which case $X^N_2(0)=N-X^N_1(0) = O(N)$.   We will be interested in letting $N\to \infty$.  In this setting, $A$ is a discrete species, $\alpha_1=0$, and $B$ is a continuous one, $\alpha_2 = 1$. Moreover, for each reaction $y_r\to y_r'$ we have $\beta_r=1$, and Assumption \ref{ass:discrete=fast}  is fulfilled. Furthermore, the limiting rate functions defined in \eqref{eq:lambda_convergence} are given by
 \begin{align*}
 	\lambda_{A+B\to2B}(x) &= \lim_{N\to \infty} N^{-1} \lambda^N_{A+B \to 2B}([N^\alpha x]) = \lim_{N\to \infty} N^{-1}\kappa_1 x_1 \floor{Nx_2} =  \kappa_1 x_1x_2 \\
 	\lambda_{B\to A}(x) &= \lim_{N\to \infty} N^{-1} \lambda^N_{B \to A}([N^\alpha x]) = \lim_{N\to \infty} N^{-1}\kappa_2  \floor{Nx_2} =  \kappa_2x_2.
\end{align*}
 
 Turning to Example \ref{ex:non_mass action} we suppose  the conserved quantity satisfies
 \[
 X^N_1(0)+X^N_2(0)+X^N_3(0)=N,
\]
where $N$ is large,  and that $X^N(0)$ is not far from the equilibrium $q^N$, which satisfies
 $$q^N_1=\frac{\kappa_0\kappa_1\kappa_3}{\kappa_2+\kappa_3}\cdot\frac{q^N_2+1}{q^N_2-1},\quad q^N_3=\frac{\kappa_1}{\kappa_2+\kappa_3}q^N_2\quad\text{and}\quad q^N_1+q^N_2+q^N_3=N.$$
 Therefore, for $N$  large $X^N_1(0)$ will be near  the value $q=\frac{\kappa_0\kappa_1\kappa_3}{\kappa_2+\kappa_3}$, while $X^N_2(0)$ and $X^N_3(0)$ go to infinity as $N\to \infty$. In this context, $A$ is a discrete species, $\alpha_1=0$, and $B$ and $C$ are continuous species, $\alpha_2 = \alpha_3 = 1$. Furthermore, $\beta_r = 1$ for each $r$, and Assumption \ref{ass:discrete=fast}  holds. In this case, the limiting rate functions defined in \eqref{eq:lambda_convergence} are given by
 \begin{align*}
 \lambda_{A+2B\to3B}(x)&=\lim_{N\to \infty} N^{-1} \kappa_0 \frac{x_1 N x_2(Nx_2 - 1)}{1+Nx_2} = \kappa_0x_1x_2,
 \end{align*}
 and similarly,
 \begin{align*}
 \lambda_{B\to C}(x)&=\kappa_1x_2,\quad\lambda_{C\to B}(x)=\kappa_2x_3,\quad\text{and}\quad\lambda_{C\to A}(x)=\kappa_3x_3.
\end{align*}

Returning to the general setting, let $\fastR$ be the set of reactions whose source complex contain a continuous species, i.e.
 \begin{equation}\label{eq:fastR}
  \fastR= \{y\to y'\in\R\,:\,\pi_c(y)\neq0\}.
 \end{equation}
 Due to \eqref{eq:lambda_convergence}, the reactions in $\fastR$ have much higher rates than the other reactions, when $N$ is large. Therefore, they give the major contribution to the dynamics of the stochastic system, and we focus on them.
 
 We define two reduced systems, one being a projection onto the discrete species space and the other being the projection onto the continuous species space of the dynamics induced by the reactions in $\fastR$. We begin by considering the projection onto the discrete species. Define
 \begin{align*}
 \pi_d(\C)&=\{\pi_d(y)\,:\, y\in\C\},\\
 \fastR_d&=\{\pi_d(y)\to\pi_d(y')\,:\, y\to y'\in\fastR \text{ and } \pi_d(y) \ne \pi_d(y')\},
 \end{align*}
 and let $\mathcal{G}_d = (\Sp_d,\pi_d(\C),\fastR_d)$ be the reaction network associated with the discrete species. For example, for both Example \ref{ex:simple} and Example \ref{ex:non_mass action}, the network associated with the discrete species is
 \[
 	A \ce{<=>} 0,
 \]
 which has a Poisson stationary distribution.
 Let $z_k \to z_k' \in \fastR_d$. For any vector $w\in\RR^{|\Sp_c|}_{\geq0}$ we define the function $\lambda_{d,k}^w : \RR^{|\Sp_d|}\to \RR_{\ge 0}$ via
 \begin{align}
 \label{eq:6789}
 \lambda^w_{d,k}(v)&=\sum_{y_r\to y'_r\in\fastR^{d,k}}\lambda_{r}(v,w),
 \end{align}
 where $\fastR^{d,k} = \{ y_r\to y'_r\in \fastR\, : \, \pi_d(y_r) =z_k\text{ and }\pi_d(y_r')=z_k'\}$. Let $\K^w_d$ be the kinetics defined by \eqref{eq:6789}, and define $\Sy^w_d=(\mathcal{G}_d,\K^w_d)$. 
 Note that the functions $\lambda_{r}$ in \eqref{eq:6789} are the \emph{limit} rate functions in \eqref{eq:lambda_convergence}. The sum in \eqref{eq:6789} is needed, as the cardinality of $\fastR^{d,k}$ is not necessarily 1. Consider for example the following modification of Example \ref{ex:simple}:
 \begin{align*}
    A + B \ce{->[\kappa_1]} 2B \ce{<=>[\kappa_3][\kappa_4]} A+2B, \qquad 
    B \ce{->[\kappa_2]} A.   
 \end{align*}
 The reactions $2B\to A+2B$ and $B\to A$ collapse to the same reaction in $\fastR_d$, and the same happens to $A+2B\to 2B$ and $A+B\to 2B$. In this case, if $w$ denotes the concentration of the species $B$, the system $\Sy^w_d$ is given by
 $$A\ce{<=>[\kappa_1w+\kappa_4w^2][\kappa_2w+\kappa_3w^2]}0.$$
 
 We make a key structural assumption on our models.
 \begin{assumption}\label{ass:complex_balanced}
 For any  $y_r\to y'_r\in \fastR$, 
 \[
 	\lambda_r(v,w)=\kappa_r(w) \frac{v!}{(v-\pi_d(y_r))!},
\]
for some functions $\kappa_r$ such that $\kappa_r(w)>0$ whenever $w>0$. It follows that for any $w>0$ the system $\Sy^w_d$ is endowed with stochastic mass action kinetics. Further, we assume that for any $w>0$ the system is complex balanced. We also require that $\G_d$ endowed with mass action kinetics is irreducible.
 \end{assumption} 

 For the systems in Examples \ref{ex:simple} and \ref{ex:non_mass action}, $\Sy^w_d$ is given, respectively, by the stochastic mass action systems
  \begin{align}
  \label{eq:67896978}
  	A\ce{<=>[\kappa_1w][\kappa_2w]}0\quad\text{and}\quad A\ce{<=>[\kappa_0w_1][\kappa_3w_2]}0,
	\end{align}
  where in the first system $w$ represents the amount of species $B$ and in the second system $w_1$ and $w_2$ represent the amounts of species $B$ and $C$, respectively. In both cases, Assumption \ref{ass:complex_balanced} holds, due to Theorem \ref{thm:deficiency_zero_iff}.

 Thanks to Assumption \ref{ass:complex_balanced} and Theorem \ref{thm:complex_balanced}, we know that for any positive $w$ there is precisely one complex balanced equilibrium of the system $\Sy^w_d$, which we denote by $q^w_d$. For example, for the first system in \eqref{eq:67896978} associated with Example \ref{ex:simple} we have $q^w_d = \kappa_2/\kappa_1$, whereas for the second system in \eqref{eq:67896978} associated with Example \ref{ex:non_mass action} we have $q^w_d = \kappa_3w_2/\kappa_0 w_1$.
 
 Assumption \ref{ass:complex_balanced} is the last structural assumption we require  for our main results. Now we impose some conditions ensuring that the systems are ``well-behaved.'' Specifically, we want that the magnitude of the counts of each species is maintened constant in a compact interval of time. In particular, we want to rule out the possibility of a blow up or of a zeroing of the concentrations of the continuous species, and we want the counts of the discrete species to be bounded, in some sense.
 
 We start by considering the projection onto the continuous species. Let
 \begin{align*}
 \pi_c(\C)&=\{\pi_c(y)\,:\, y\in\C\},\\
 \fastR_c&=\{\pi_c(y)\to\pi_c(y')\,:\, y\to y'\in\fastR \text{ and } \pi_c(y) \ne \pi_c(y')\},
 \end{align*}
 and define $\mathcal{G}_c = (\Sp_c,\pi_c(\C),\fastR_c)$. We now define the kinetics for the network $\mathcal{G}_c$.  
 For any $y_r\to y'_r\in \fastR$, define 
 \begin{equation}\label{eq:fast_kinetics}
  \widetilde{\lambda}_{r}(v,w)=\kappa_r(w)v^{\pi_d(y_r)},
 \end{equation}
 where the $\kappa_r(w)$ are as in Assumption \ref{ass:complex_balanced}.  Note  we are assuming that the intensities $\widetilde{\lambda}_{r}(v,w)$ take the form of deterministic mass action kinetics and not stochastic mass action. 
For any $z_k\to z_k'\in \fastR_c$, we then define the function $\lambda_{c,k}:\RR^{|\Sp_d|}\to \RR_{\ge 0}$ via
 \begin{equation}\label{eq:definition_rates_continuous_system}
  \lambda_{c,k}(w)=\sum_{y_r\to y'_r\in\fastR^{c,k}}\mathbbm{1}_{\{w>0\}}\widetilde{\lambda}_{r}(q^w_d,w),
 \end{equation}
 where $\fastR^{c,k} = \{ y_r\to y'_r\in \fastR\, : \, \pi_c(y_r) =z_k\text{ and }\pi_c(y_r')=z_k'\}$. Let $\K_c$ denote the kinetics defined by the above rate functions, and define $\Sy_c=(\mathcal{G}_c,\K_c)$. Finally, fix a finite time $T>0$ and a point $X_0\in\RR^{|\Sp|}_{>0}$.
 
 \begin{assumption}\label{ass:existence_limit}
 Assume that the deterministic solution $z(t)$ of the system $\Sy_c$, with initial condition $\pi_c(X_0)$, exists for any $t\in[0,T]$. Moreover, assume that for all $t\in[0,T]$ we have $z(t)>0$.
 \end{assumption}

 \begin{remark}
  A slightly more general scenario than that given by Assumption \ref{ass:existence_limit} could be considered.  Specifically, we could allow $z(t)$ to be equal to zero for some $t\in[0,T]$.  In fact, for our purposes it is enough that (i)  the deterministic solution $z(t)$ exists and (ii) that the reaction rates $\lambda_r$ are of the form described in Assumption \ref{ass:complex_balanced} and (iii) $\Sy^w_d$ is complex balanced for any $w$ in a neighborhood of $\{z(t)\}_{t\in[0,T]}$ (relative to the non-negative orthant). For Corollary \ref{cor:discrete_species_unary} to hold, we will further need to assume that for $w$ varying in such a neighborhood, the quantities $\kappa_r(w)$ introduced in Assumption \ref{ass:complex_balanced} are bounded from below by a positive constant.
 \end{remark}

  Consider Example \ref{ex:simple}. In this case, for any $w\in\RR_{>0}$ we have $q^w_d=\kappa_2/\kappa_1$, and the system $\Sy_c$ is given by
  \begin{equation}\label{eq:ex_simple_continuous_system}
  0\ce{<-[\kappa_2]}B\ce{->[\kappa_2]}2B.
  \end{equation}
  Hence, in this case the deterministic solution $z(t)$ is constantly equal to $\pi_c(X_0)$.
  
  Consider now Example \ref{ex:non_mass action}. Here, for any $w\in\RR^2_{>0}$, $q^w_d$ is given by $\kappa_3 w_2/\kappa_0 w_1$. Therefore, the system $\Sy_c$ is
  \begin{gather*}
   2B \xrightarrow{\phantom{\kappa_3}} 3B,\qquad 
   B\ce{<=>[\kappa_1][\kappa_2]} C \xrightarrow{\kappa_3} 0,
  \end{gather*}
  with
  $\lambda_{2B\to3B}(w)=\kappa_3w_2.$
  When deterministically modeled, the dynamics of the system is equivalent to that of the deterministic mass action system
  \begin{equation}\label{eq:ex_not_mass_action_continuous_system}
   B\ce{<=>[\kappa_1][\kappa_2+\kappa_3]} C
  \end{equation}
  and it can be easily shown that Assumption \ref{ass:existence_limit} holds, since $\pi_c(X_0)>0$.
 
 \begin{assumption}\label{ass:moments_bound}
 There exists a locally bounded function $\psi\colon\RR^{|\Sp_d|}_{\geq0}\to\RR_{\geq1}$ satisfying\vspace*{.5em}\\ 
 \begin{enumerate*}[label=(\roman*), itemjoin={{\quad\text{and}\quad}}, before={{\phantom{.}\hfill}}, after={{,\hfill\phantom{.}}}]
  \item\label{eq:moments_bound_lim_to_infty} $\displaystyle\lim_{\|v\|\to\infty}\psi(v)=\infty$
  \item\label{eq:moments_bound_sup_bounded}  $\displaystyle\sup_{N\in\NN_{>0}} \sup_{t\in[0,T]}E\prq{\psi(\Xd^N(s))}<\infty$
 \end{enumerate*}\vspace*{.5em}\\
 such that for any reaction $y_r\to y_r'\in\R$ and any compact set $\Gamma\subset\RR^{|\Sp_c|}_{\geq0}$\vspace*{.5em}\\
 \begin{enumerate*}[label=(\roman*), before={{\phantom{.}\hfill}}, after={{.\hfill\phantom{.}}}] 
 \setcounter{enumi}{2}
  \item\label{eq:moments_bound_limsup_to_zero} $\displaystyle\sup_{w\in\Gamma}\sup_{v\in\RR_{>0}^{|\Sp_d|}}\frac{\lambda_r(v,w)}{\psi(v)}<\infty$
 \end{enumerate*}\vspace*{.5em}\\
 \end{assumption}
 \begin{remark}\label{rem:moments_bound_mass_action}
  For any reaction $y_r \to y_r'\in \mathcal{R}$ and any $N \in \mathbb{N}$ assume that $\lambda_r$
  has at most polynomial growth rate in $v$. Let $p\colon\RR^{|\Sp_d|}\to\RR$ be a polynomial satisfying
  $$\max_{y_r\to y_r'\in\R}\limsup_{\|v\|\to\infty}\frac{\lambda_r(v,w)}{p(v)}=0\quad\text{for any }w\in\RR_{\geq0}^{|\Sp_c|},$$
  and let $d$ be the degree of $p$. Then, a candidate for $\psi$ is
  $$\psi(v)=1+\sum_{S_i\in\Sp_d}v_i^d.$$
  With this choice, $\psi$ automatically satisfies \ref{eq:moments_bound_lim_to_infty} and \ref{eq:moments_bound_limsup_to_zero}.
 \end{remark}

 We now state and prove the theorem which provides the backbone for our results in the setting of absolute concentration robustness.
 
 \begin{theorem}\label{thm:weak_convergence}
 If Assumptions \ref{ass:discrete=fast} to \ref{ass:moments_bound} hold, then 
 $$\sup_{t\in[0,T]}\abs{N^{-1}\Xc^N(t)-z(t)}\xrightarrow[N\to\infty]{}0$$
 in probability, where $z$ is defined as in Assumption \ref{ass:existence_limit}. 
 
 Moreover, if $\Pois(q)$ denotes a product-form Poisson distribution with parameter $q$, then for any continuous function $g\colon\RR^{|\Sp_d|}_{\geq0}\to\RR$ satisfying
 \begin{equation}\label{eq:g_bounded}
  \limsup_{\|v\|\to\infty}\frac{|g(v)|}{\psi(v)}=0
 \end{equation}
 we have
 $$\sup_{t\in[0,T]}\abs{\int_0^{t}\Big(g(\Xd^N(s))-E[g(J^{z(s)})]\Big)ds}\xrightarrow[N\to\infty]{}0$$
 in probability and for each $\delta \in(0,T)$
\begin{equation}\label{eq:fixed_time_conv}
\sup_{t\in[\delta,T]}\abs{E[g(\Xd^N(t))]-E[g(J^{z(t)})]}\xrightarrow[N\to\infty]{}0.
\end{equation}
 where $J^{z(s)}\sim\Pois(q_d^{z(s)})$.
 \end{theorem}
 \begin{remark}\label{rem:general_stationary_distribution}
  If Assumption \ref{ass:complex_balanced} is not satisfied, but we know that for any $w\in\RR_{>0}^{|\Sp_c|}$ the stochastic system $\Sy^w_d$ possesses a unique stationary distribution $\mu^w$ with
  $$E_{\mu^w}[\lambda_r(v,w)]=\sum_{v\in\NN^{\Sp_d}}\lambda_r(v,w)\mu^w(v)<\infty\quad\text{for every }y_r\to y_r'\in\R,$$
  then Theorem \ref{thm:weak_convergence} still holds, provided that every occurrence of $\lambda_r(q_d^w,w)$ is replaced by $E_{\mu^w}[\lambda_r(v,w)]$. The proof of Theorem \ref{thm:weak_convergence}, with small changes, also covers this generalization. In regard to this broader setting, see also the results in \cite{kurtz:rescale,popovic:rescale}.
 \end{remark}
 \begin{proof}
 For the sake of simplicity, throughout the proof we will write $t$ instead of $t\wedge T$, but it is always implicitly assumed that $t\in[0,T]$.
 
 We follow the arguments of \cite{kurtz:rescale,popovic:rescale}, which rely on the techniques developed in \cite{kurtz:averaging}. We will first prove the theorem under the assumption that
 \begin{equation}\label{eq:stochastically_bounded}
  \sup_N P\pr{\sup_{t\in[0,T]}N^{-1}\|\Xc^N(t)\|_\infty>M}=0,
 \end{equation}
 for a certain constant $M$ satisfying
 $$\sup_{t\in[0,T]}\|z(t)\|_\infty<M.$$
 The above assumption holds when a positive linear combination of the species is conserved, or if we study the process up to the time when the concentration of a continuous species exceeds a given threshold.  We will then drop the assumption \eqref{eq:stochastically_bounded}.
 
 Define the \emph{occupation measures} $\Gamma^N$ on $\RR^{|\Sp_d|}\times[0,T]$ by
 $$\Gamma^N(D\times [a,b])=\int_a^b \mathbbm{1}_D(\Xd^N(s))ds.$$
 Note that
\begin{equation}\label{eq:gamma_N}
 d\Gamma^N(v,s)=d\gamma^N_s(v)ds,
 \end{equation}
 where $\gamma^N_s=\delta_{\Xd^N(s)}$, with $\delta_x$ denoting the usual Dirac measure on $\RR^{|\Sp_d|}$.
 By part \ref{eq:moments_bound_sup_bounded} in Assumption \ref{ass:moments_bound}, we have that for any $\varepsilon>0$ there exists a constant $M_\varepsilon>0$ such that
 \[
 \sup_{N\in\NN_{>0}} \sup_{s\in[0,T]}P(\psi(\Xd^N(s)) <M_\varepsilon)>1-\varepsilon.
 \]
 It follows that
 $$E\prq{\Gamma^N(\psi^{-1}([0,M_\varepsilon])\times [0,T])}\geq (1-\varepsilon)T.$$
 By part \ref{eq:moments_bound_lim_to_infty} in Assumption \ref{ass:moments_bound}, we have that $\psi^{-1}([0,M_{\varepsilon}])$ is compact, and Lemma 1.3 in \cite{kurtz:averaging} implies that the sequence of random measures $\Gamma^N$ is relatively compact (see also Lemma 2.9 in \cite{kurtz:central}). Let $\Gamma$ be a weak limit point.

 Consider the generator $L^N$ for the process $N^{-\alpha}X^N(t)$, defined by
 $$L^Nf(x)=\sum_{y_r\to y_r'\in\R}\lambda^N_r(N^{\alpha}x)\Big(f\pr{x+N^{-\alpha}(y'-y)}-f(x)\Big),\quad\text{for }N^{\alpha}x\in\NN^{|\Sp|}.$$
 From the generator $L^N$ we can obtain two generators, one related to the limiting behaviour of the concentrations of the continuous species (whose changes take place at the time scale $t$) and the other one related to the discrete species (whose changes take place at the time scale $N^{-1}t$).
 For any function $h\in C^2_c(\RR^{|\Sp_c|})$ and $x\in\NN^{|\Sp_d|}\times\RR^{|\Sp_c|}$, define
 \begin{align*}
 L_ch(x)=&\lim_{N\to\infty}L^N(h\circ\pi_c)(N^{-\alpha}[N^{\alpha}x])\\
 =&\sum_{y_r\to y_r'\in\R}\lim_{N\to\infty}\lambda^N_r([N^{\alpha}x])\Big(h\Big(\pi_c(N^{-\alpha}[N^{\alpha}x])+N^{-1}\pi_c(y'-y)\Big)-h\Big(\pi_c(N^{-\alpha}[N^{\alpha}x])\Big)\Big)\\
 =&\sum_{y_r\to y_r'\in\fastR}\lambda_r(x)\pi_c(y'-y)\cdot \nabla h(\pi_c(x)),
 \end{align*}
 where we made use of Assumption \ref{ass:discrete=fast} to compute the limit, and $\cdot$ denotes the scalar product. Note that $L_ch\in C^1_c(\RR^{|\Sp|})$ and that $L^N(h\circ\pi_c)(N^{-\alpha}[N^{\alpha}x])$ converges uniformly to $L_ch(x)$ on $x$ (which follows from the fact that $h$ has compact support). $L_c$ can be interpreted as the generator of the limiting behaviour of the concentrations of the continuous species.
 
 On the other hand, for any function $g\in C_c(\NN^{|\Sp_d|})$ and $x\in\NN^{|\Sp_d|}\times\RR^{|\Sp_c|}$, define
 \begin{align*}
 L_dg(x)=&\lim_{N\to\infty}N^{-1}L^N(g\circ\pi_d)(N^{-\alpha}[N^{\alpha}x])\\
 =&\sum_{y_r\to y_r'\in\R}\lim_{N\to\infty}N^{-1}\lambda^N_r([N^{\alpha}x])\Big(g\Big(\pi_d(x)+\pi_d(y'-y)\Big)-g\Big(\pi_d(x)\Big)\Big)\\
 =&\sum_{y_r\to y_r'\in\fastR}\lambda_r(x)\Big(g\Big(\pi_d(x)+\pi_d(y'-y)\Big)-g\Big(\pi_d(x)\Big)\Big).
 \end{align*}
 The convergence is uniform in $x$. $L_d$ can be interpreted as the generator of the limiting behavior of the discrete species on the timescale $N^{-1} t$. For any $w\in\RR^{|\Sp_c|}_{\geq0}$ we can define the operator $L^w_d$ by
 $$L^w_dg(v)=L_dg(v,w)\quad\forall v\in\NN^{|\Sp_d|},$$
 which corresponds to the generator of the system $\Sy^w_d$.
 
 For any $h\in C^2_c(\RR^{|\Sp_c|})$, the process
 \begin{equation}\label{eq:M^N_h}
  M^N_h(t)=h(N^{-1}\Xc^N(t))-h(N^{-1}\Xc^N(0))-\int_0^tL^N(h\circ\pi_c)(N^{-\alpha}X^N(s))ds
 \end{equation}
 is a martingale. Let
 $$\delta^N_h(t)=\int_0^t\pr{L_ch-L^N(h\circ\pi_c)}(N^{-\alpha}X^N(s))ds.$$
 By the uniform convergence of $L^N(h\circ\pi_c)(N^{-\alpha}[N^{\alpha}x])$ to $L_ch(x)$ with respect to $x$, we have that
 \begin{equation}\label{eq:lim_delta}
  \lim_{N\to\infty}E\prq{\sup_{t\in[0,T]}|\delta^N_h(t)|}=0.
 \end{equation}
 Moreover, \ref{eq:moments_bound_sup_bounded} and \ref{eq:moments_bound_limsup_to_zero} in Assumption \ref{ass:moments_bound}, together with \eqref{eq:stochastically_bounded}, imply that
 \begin{equation}\label{eq:L_ch_stochastic_bound}
  \sup_N E\prq{\int_0^T|L_ch(N^{-\alpha}X^N(s))|ds}<\infty.  
 \end{equation}
 Since \eqref{eq:M^N_h} is a martingale, tightness of the processes $\Xc^N$ stopped at $T$ follows from \eqref{eq:stochastically_bounded}, \eqref{eq:lim_delta} and \eqref{eq:L_ch_stochastic_bound}, and by Theorems 3.9.1 and 3.9.4 in \cite{ethier_kurtz}.\footnote{To be specific, the processes to which we apply Theorem 3.9.4 in \cite{ethier_kurtz}  are,
 \begin{align*}
  Z_n(t)=L_ch(N^{-\alpha}X^N(t)), \quad \text{and} \quad
  Y_n(t)=&M^N_h(t)+\int_0^tZ_n(s)ds,
 \end{align*}
 where $Z_n,Y_n$ is the notation utilized in \cite{ethier_kurtz}.
 This notation is now dropped throughout the remainder of this paper.}
 
 Let $(W, \Gamma)$ be a weak limit of $(\Xc^N(\cdot\wedge T),\Gamma^N)$. By the same arguments as in the proof of Theorem 2.1 in \cite{kurtz:averaging}, we have that
 $$M_h(t)=h(W(t))-h(W(0))-\int_{\RR^{|\Sp_d|}\times[0,t]}L_ch(v,W(s))d\Gamma(v,s)$$
 is a martingale. On the other hand, for any $g\in C^2_c(\RR^{|\Sp_d|})$, the process
 $$\widehat{M}^{N}_g(t)=N^{-1}\prq{g(\Xd^N(t))-g(\Xd^N(0))-\int_0^tL^N(g\circ\pi_d)(N^{-\alpha}X^N(s))ds}$$
 is also a martingale. Since the function $g$ is bounded and by using uniform convergence of $N^{-1}L^N(g\circ\pi_d)$ to $L_dg$, we have that a weak limit point for $\widehat{M}^{N}_g$ is given by
 $$\widehat{M}_g(t)=-\int_{\RR^{|\Sp_d|}\times[0,t]}L_dg(v,W(s))d\Gamma(v,s),$$
 which is therefore a martingale. By \eqref{eq:gamma_N} we have
 \begin{equation}\label{eq:gamma}
  d\Gamma(v,s)=d\gamma_s(v)ds,
 \end{equation}
 for a family of measures $\gamma_s$. Therefore, $\widehat{M}_g(t)$ is continuous and for any $t_1<t_2$
 $$\int_{\RR^{|\Sp_d|}\times[t_1,t_2]}\left|L_dg(v,W(s))\right|d\Gamma(v,s)\leq (t_2-t_1)\sup_{x\in\RR^{\Sp}}\left|L_dg(x)\right|,$$
 which implies that $\widehat{M}_g(t)$ has finite variation paths. This in turn implies that $\widehat{M}_g(t)$ is constantly equal to $\widehat{M}_g(0)=0$ for any $t\in[0,T]$ with probability one. Therefore, almost surely, for almost every $s\in[0,T]$
 $$\int_{\RR^{|\Sp_d|}}L^{W(s)}_dg(v)d\gamma_s(v)=0,$$
 where $\gamma_s$ is as in \eqref{eq:gamma}.
 Since $C^2_c(\RR^{|\Sp|})$ is separable, we have that, for almost every $s\in[0,T]$,
 \[
 	\int_{\RR^{|\Sp_d|}}L^{W(s)}_dg(v)d\gamma_s(v)=0,\quad\forall g\in C^2_c(\RR^{|\Sp_d|}).
\]
 Thus, for almost every $s$ such that $W(s)>0$ with probability one, the measure $\gamma_{s}$ is equal to the unique stationary distribution of the system $\Sy^{W(s)}_d$. Due to Assumption \ref{ass:complex_balanced} and by Theorem 4.1 in \cite{anderson:product-form}, it corresponds to the product form Poisson distribution $\Pois(q^{W(s)}_d)$. Therefore we can write
 $$\gamma_s=\mathbbm{1}_{\{W(s)>0\}}\Pois(q^{W(s)}_d)+(1-\mathbbm{1}_{\{W(s)>0\}})\gamma_s.$$
 The weak limit process $W(t)$ is a solution to the martingale problem
 \begin{align*}
 M_h(t)=&h(W(t))-h(W(0))-\int_0^t\int_{\RR^{|\Sp_d|}}L_ch(v,W(s))d\gamma_s(v)ds\\
 =&h(W(t))-h(W(0))-\sum_{y_r\to y_r'\in\fastR}\int_0^t\mathbbm{1}_{\{W(s)>0\}}\lambda_r(q^{W(s)}_d,W(s))\pi_c(y'-y)\cdot \nabla h(W(s))ds\\
 &-\int_0^t(1-\mathbbm{1}_{\{W(s)>0\}})\int_{\RR^{|\Sp_d|}}L_ch(v,W(s))d\gamma_s(v)ds,
 \end{align*}
 where in the last equality we used that if $V$ is a product form Poisson random variable with mean $q$, then for any non-negative integer vector $y$ of the same dimension of $V$ 
 $$E\prq{\frac{V!}{(V-y)!}}=q^y.$$
 This in turn implies that $E\prq{\lambda_r(V,w)}=\widetilde{\lambda}_r(q,w)$, where $\widetilde{\lambda}_r$ are as defined in \eqref{eq:fast_kinetics}. By Assumption \ref{ass:existence_limit}, $W(t)$ is uniquely determined by the solution to
 \begin{align*}
 W(0)=&\pi_c(X_0)\\
 W(t)=&W(0)+\sum_{y_r\to y_r'\in\fastR}\int_0^t\mathbbm{1}_{\{W(s)>0\}}\widetilde{\lambda}_r(q^{W(s)}_d,W(s))\pi_c(y'-y)ds+\\
 &+\int_0^t(1-\mathbbm{1}_{\{W(s)>0\}})\int_{\RR^{|\Sp_d|}}L_cid(v,W(s))d\gamma_s(v)ds\\
 =&W(0)+\sum_{y_r\to y_r'\in\fastR}\int_0^t\widetilde{\lambda}_r(q^{W(s)}_d,W(s))\pi_c(y'-y)ds,
 \end{align*}
 which is given by $z(t)$. 
 The first part of the theorem is therefore proved. The second part follows from Lemma 2.9 in \cite{kurtz:central}.
 
 To prove that the first two parts of the theorem hold without assuming \eqref{eq:stochastically_bounded}, fix two positive constants $M$ and $\delta$ with the property
 $$\sup_{t\in[0,T]}\|z(t)\|_\infty<M-\delta,$$
 and consider the stopping time
 $$\tau^N=\inf\{t\in[0,T]\,:\,N^{-1}\|\Xc^N(t)\|_\infty>M\}.$$
 Note that $N^{-1}\Xc(\tau^N)$ is also uniformly bounded in $N$, as
 $$N^{-1}\|\Xc^N(\tau^N)\|_\infty\leq M+\max_r\{\|\xi_r\|_\infty\}=M'.$$
 Therefore, \eqref{eq:stochastically_bounded} holds up to time $\tau^N$ with the constant $M'$, which means that for any $\varepsilon$ there exists $N_\varepsilon$ such that for any $N>N_\varepsilon$
 $$P\pr{\sup_{t\in[0,\tau^N]}\|N^{-1}\Xc^N(t)-z(t)\|>\delta}<\varepsilon.$$
 Since $\{\tau^N<T\} \subset \{\|N^{-1}\Xc^N(\tau^N)-z(\tau^N)\|>\delta\}\subset \{\sup_{t\in[0,\tau^N]}\|N^{-1}\Xc^N(t)-z(t)\|>\delta\}$, it follows that for $N$ large enough $P(\tau^N<T)<\varepsilon$. 
 We conclude our argument by the arbitrariness of $\varepsilon$, since for $N$ large enough
 \begin{align*}
  P\pr{\sup_{t\in[0,T]}\|N^{-1}\Xc^N(t)-z(t)\|>\eta}&\leq P\pr{\sup_{t\in[0,\tau^N]}\|N^{-1}\Xc^N(t)-z(t)\|>\eta}+\varepsilon,\\
  P\pr{\sup_{t\in[0,T]}\int_0^{t}\Big(g(\Xd^N(s))-E[g(J^{z(s)})]\Big)ds>\eta} &\leq P\pr{\sup_{t\in[0,\tau^N]}\int_0^{t}\Big(g(\Xd^N(s))-E[g(J^{z(s)})]\Big)ds>\eta}+\varepsilon.
 \end{align*}
 
 To show the last part of the theorem, namely \eqref{eq:fixed_time_conv}, consider a converging sequence $\{t_N\}$ in $[\delta,T]$ with limit $t$. Fix a positive real number $r_0$ and for any $r\leq r_0$ and any $N>r_0/\delta$ define
 $$Y_{r_0}^N(r)=\Xd^N\pr{t_N+\frac {r-r_0}{N}}.$$
 By Assumption \ref{ass:moments_bound} \ref{eq:moments_bound_sup_bounded}, the sequence $Y_{r_0}^N(0)$ is relatively compact. Let $Y_{r_0}(r)$ be a process with generator $L^{z(t)}_d$ and such that $Y_{r_0}(0)$ is the weak limit point of a subsequence $Y_{r_0}^{N_m}(0)$. By using again Assumption \ref{ass:moments_bound} \ref{eq:moments_bound_sup_bounded} and the weak convergence of $N^{-1}\Xc^N(s)$ to $z(s)$,  for any function $g:\NN^{\Sp_d}\mapsto \RR_{\geq0}$ that satisfies \eqref{eq:g_bounded} we have 
 \begin{align*}
  \lim_{m\to\infty}E\prq{g\pr{\Xd^{N_m}\pr{t_{N_m}}}}=&\lim_{m\to\infty}E\prq{g\pr{Y_{r_0}^{N_m}(r_0)}}\\
  =&\lim_{m\to\infty}E\prq{g\pr{Y_{r_0}^{N_m}(0)}}\\
  &+\lim_{m\to\infty}\int_0^{r_0}{N_m}^{-1}E\prq{L^{N_m}(g\circ\pi_d)\pr{Y_{r_0}^{N_m}(r),\Xc^{N_m}\pr{t_{N_m}+\frac {r-r_0}{{N_m}}}}}dr\\
  =&E\prq{g\pr{Y_{r_0}(0)}}+\int_0^{r_0}E\prq{L^{z(t)}_dg(Y_{r_0}(r))}dr.
 \end{align*}
 Note that the left-hand side of the above equation does not depend on $r_0$. Moreover, Assumption \ref{ass:moments_bound}\ref{eq:moments_bound_sup_bounded} implies that the limit $E\prq{g\pr{Y_{r_0}(0)}}$ is bounded from above, independently of $r_0$. Therefore, we can consider a weakly convergent sequence $g\pr{Y_{r^n_0}(0)}$ with $r^n_0\to\infty$ as $n\to\infty$. Letting $g(Y(0))$ denote the weak limit, we have
 $$\lim_{m\to\infty}E\prq{g\pr{\Xd^{N_m}\pr{t_{N_m}}}}=E\prq{g(Y(0))}+\int_0^\infty E\prq{L^{z(t)}_dg(Y(r))}dr.$$
 Since $L^{z(t)}$ is the generator of a stochastic process with stationary distribution $\Pois(q_d^{z(t)})$, we must have $$\lim_{m\to\infty}E\prq{g\pr{\Xd^{N_m}\pr{t_{N_m}}}}=E\prq{g(J^{z(t)}_q)}$$
 Moreover, since the limit does not depend either on the particular subsequence $N_m$ or on the converging sequence $t_N$, we conclude
 $$\sup_{t\in[\delta,T]}\abs{E[g(\Xd^N(t))]-E[g(J^{z(t)})]}\xrightarrow[N\to\infty]{}0.$$
 which is \eqref{eq:fixed_time_conv}.
 \end{proof}
 
 In some cases Assumption \ref{ass:moments_bound} can be difficult to check, even if it seems natural for the analysed system. For this reason, we state here a corollary of Theorem \ref{thm:weak_convergence} concerning a particular case for which Assumption \ref{ass:moments_bound} is automatically satisfied.  
 
 \begin{corollary}\label{cor:discrete_species_unary}
  Assume Assumptions \ref{ass:discrete=fast} to \ref{ass:existence_limit} hold. Assume also that
  \begin{equation}\label{eq:extra_condition_main_corollary_above}
   N^{-\beta_r}\lambda_r^N(v,[Nw])\leq \widehat{h}_r(w)v^{\pi_d(y_r)}\quad\text{for any }y_r\to y_r'\in\R,
  \end{equation}
  for some continuous positive functions $\widehat{h}_r\colon\RR^{|\Sp_c|}\to\RR$, and that
  \begin{equation}\label{eq:extra_condition_main_corollary_below}
   N^{-\beta_r}\lambda_r^N(v,[Nw])\geq h_r(w)w^{\pi_c(y_r)}v^{\pi_d(y_r)}\quad\text{for any }y_r\to y_r'\in\fastR,
  \end{equation}  
  for some continuous positive functions $h_r\colon\RR^{|\Sp_c|}\to\RR$. Furthermore, assume that in the support of any complex $y\in\C$ at most one discrete species appears, and its stoichiometric coefficient is 1. Then, for any continuous function $g\colon\RR^{|\Sp_d|}_{\geq0}\to\RR$ with at most polynomial growth rate we have 
  $$\sup_{t\in [0,T]}\abs{\int_0^t\Big(g(\Xd^N(s))-E[g(J^{z(s)})]\Big)ds}\xrightarrow[N\to\infty]{}0$$
  in probability, where $J^w\sim\Pois(q^w_d)$ and $z(t)$ is as in Assumption \ref{ass:existence_limit}. Moreover, for any continuous bounded function $\varphi\colon\RR^{|\Sp_d|}_{\geq0}\to\RR$ and any $\delta\in(0,T)$, we have
  $$\sup_{t\in[\delta,T]}\abs{E[\varphi(\Xd^N(t))]-E[\varphi(J^{z(t)})]}\xrightarrow[N\to\infty]{}0.$$
 \end{corollary}
  \begin{remark}\label{rem:mass_action}
  Assume that $\K^N$ is mass action kinetics and at most one discrete species appear in the support of any complex $y\in\C$. If the rate constants are rescaled according to \eqref{eq:rate_constants_convergence},  then\eqref{eq:extra_condition_main_corollary_above} and \eqref{eq:extra_condition_main_corollary_below} are automatically satisfied. Indeed,
  $$N^{-\beta_r}\lambda^N_r(v,[Nw])=(\kappa_r+\varepsilon^1_N)(w^{\pi_c(y_r)}+\varepsilon^2_N(w))v^{\pi_d(y_r)},$$
  for some sequence $\varepsilon^1_N$ converging to zero, and some function $\varepsilon^2_N(w)$ converging to zero uniformly on $w$. In this regard, see Remark \ref{rem:rate_constants_convergence}.
 \end{remark}
 \begin{proof}
  By Assumption \ref{ass:existence_limit}, we can choose two positive constants $m<M$ such that
  \begin{equation}\label{eq:m_and_M}
   \inf_{S_i\in\Sp_c, t\in[0,T]}z_i(t)>m\quad\text{and}\quad\sup_{S_i\in\Sp_c, t\in[0,T]}z_i(t)<M.
  \end{equation}
  For any $N\in\NN_{>0}$, consider the function $\pi^N:\RR^{|\Sp|}\to\RR^{|\Sp|}$ defined by
  $$(\pi^N(x))_i=\begin{cases}
                   (Nm\vee x_i)\wedge NM &\text{if }S_i\in\Sp_c\\
                   x_i                   &\text{otherwise.}
                  \end{cases}$$
  Consider the modified family of kinetics $\overline{\K}^N$ defined by
  $$\overline{\lambda}^N_r(x)=\lambda^N_r\pr{\pi^N(x)}\quad\text{for }x\in\NN^{|\Sp|}.$$
  We have
  $$\lim_{N\to\infty}N^{-\beta_r}\overline{\lambda}^N_r([N^\alpha x])=\overline{\lambda}_r(x)=\lambda_r\pr{\pi^1(x)},$$
  where the limit is uniform on compact sets. Furthermore, $\overline{\lambda}_r$ is locally Lipschitz. Let $\overline{\K}$ be the kinetics defined by the functions $\overline{\lambda}_r$.
  Our first aim is to prove that for a suitable choice of $\psi$, Assumption \ref{ass:moments_bound} holds for the modified kinetics. We then apply Theorem \ref{thm:weak_convergence} to the reaction systems $(\G,\overline{\K}^N)$. Let $\overline{X}^N(t)$ denote the stochastic process associated with $(\G,\overline{\K}^N)$. Define
  $$\sigma^N(t)=\sum_{S_i\in\Sp_d}\overline{X}_i^N(t)=\sum_{S_i\in\Sp_d}X_i^N(t)\quad\text{and}\quad\Delta_r=\sum_{S_i\in\Sp_d}\xi_{ri}.$$
  Since the complexes $y$ are non-negative vectors, we have
  $$\Delta_r=\sum_{S_i\in\Sp_d}(y'_{ri}-y_{ri})=\|\pi_d(y_r')\|_1-\|\pi_d(y_r)\|_1.$$
  By hypothesis, for any complex $y\in\C$ we have $\|\pi_d(y)\|_1\leq 1$, which implies that $-1\leq \Delta_r\leq 1$ for any $y_r\to y_r'\in\R$.
  Moreover, we have
  \begin{equation}\label{eq:delta_norm_L1}
   \begin{array}{rcl}
    \Delta_r= 1 &\Longrightarrow& \|\pi_d(y_r)\|_1=0\\
    \Delta_r= -1&\Longrightarrow& \|\pi_d(y_r)\|_1=1.
   \end{array}  
  \end{equation}
 Furthermore,
  \begin{align}\label{eq:sigma_first_reduction}
   \sigma^N(t)&=\sigma^N(0)+\sum_{y_r\to y_r'\in\R}\Delta_rY_r\pr{\int_0^t\overline{\lambda}^N_r\pr{\overline{X}^N(s)}ds}\notag\\
   &=\sigma^N(0)+\sum_{\substack{y_r\to y_r'\in\R\\ \Delta_r=1}}Y_r\pr{\int_0^t\overline{\lambda}^N_r\pr{\overline{X}^N(s)}ds}-\sum_{\substack{y_r\to y_r'\in\R\\ \Delta_r=-1}}Y_r\pr{\int_0^t\overline{\lambda}^N_r\pr{\overline{X}^N(s)}ds}
  \end{align}
  Define
  $$M^*=\max_{y_r\to y_r'\in\R}\max_{m\leq w\leq M}\widehat{h}_r(w)\qquad\text{and}\qquad m^*=\min_{y_r\to y_r'\in\R}\min_{m\leq w\leq M}h_r(w)w^{\pi_c(y_r)},$$
  which are both positive constants.
  By \eqref{eq:delta_norm_L1} and \eqref{eq:extra_condition_main_corollary_above}, whenever $\Delta_r= 1$ we have
  $$\overline{\lambda}^N_r(v,Nw)\leq NM^*.$$
  On the other hand, by \eqref{eq:extra_condition_main_corollary_below}, if $y_r\to y_r'\in\fastR$ and $\Delta_r=-1$ then
  $$\overline{\lambda}^N_r(v,Nw)\geq N^{\beta_r}m^*v_i\quad\text{for some }S_i\in\Sp_d.$$
  By Assumption \ref{ass:discrete=fast}, all the species in $\Sp_d$ are either produced or consumed by a reaction in $\fastR_d$. Moreover, by Assumption \ref{ass:complex_balanced} the system $\Sy^w_d$ is complex balanced, which implies that it is weakly reversible by Theorem \ref{thm:complex_balanced}. Therefore, all the species in $\Sp_d$ are both produced and consumed by some reaction in $\fastR_d$. In particular, for any $S_i\in\Sp_d$ we can choose a reaction $y_{r(i)}\to y_{r(i)}'\in\fastR$ such that $\xi_{r(i)i}=-1$. We have that
  $$\sum_{\substack{y_r\to y_r'\in\R\\ \Delta_r=-1}}Y_r\pr{\int_0^t\overline{\lambda}^N_r\pr{\overline{X}^N(s)}ds}\geq \sum_{S_i\in\Sp_d}Y_{r(i)}\pr{\int_0^t Nm^*\overline{X}_i^N(s)ds}.$$
  Then, from \eqref{eq:sigma_first_reduction} it follows that
  $$\sigma^N(t)\leq \sigma^N(0)+\sum_{y_r\to y_r'\in\R}Y_r\pr{\int_0^tNM^*ds}-\sum_{S_i\in\Sp_d}Y_{r(i)}\pr{\int_0^t Nm^*\overline{X}_i^N(s)ds},$$
  which implies that $\sigma^N(t)$ is stochastically bounded by
  $$B(t)=\sigma^N(0)+Y'\Big(|\R|\cdot M^*t\Big)-Y''\pr{\int_0^t m^*B(s)ds},$$
  where $Y'$ and $Y''$ are two i.i.d.\ unit-rate Poisson processes. That is, for all $u \ge 0$
  $$P\pr{\sup_{t\in[0,T]} \sigma^N(t) \ge u} \le P\pr{\sup_{t\in[0,T]} B(Nt) \ge u}.$$
  It follows that
  $$\sup_{t\in[0,T]}E[\sigma^N(t)^n]\leq \sup_{t\in[0,T]}E[B(Nt)^n].$$
  Then, for any $n\in\NN$
  \begin{equation}\label{eq:sigma_uniform_bound}
   \sup_{\substack{t\in[0,T]\\N\in\NN_{>0}}} E[\sigma^N(t)^n]\leq\sup_{t\in[0,\infty)}E[B(t)^n]<\infty.
  \end{equation}
  The last inequality is due to ergodic properties of the birth-death process $B(t)$. For any $n>1$, we define $\psi_n\colon\RR^{|\Sp_d|}_{\geq0}\to\RR_{\geq1}$ via
  $$\psi_n(v)=1+\sum_{S_i\in\Sp_d}v_i^n.$$
  Due to \eqref{eq:extra_condition_main_corollary_above}, $\lambda^N_r$ has at most polynomial growth rate in $v$, for any reaction $y_r\to y_r'\in\R$ and any $N\in\NN$. By \eqref{eq:lambda_convergence}, the condition \eqref{eq:extra_condition_main_corollary_above} also implies that the rate functions $\lambda_r$ have at most polynomial growth rate.
  By Remark \ref{rem:moments_bound_mass_action}, if $n$ is large enough, $\psi_n$ satisfies \ref{eq:moments_bound_lim_to_infty} and \ref{eq:moments_bound_limsup_to_zero} in Assumption \ref{ass:moments_bound}. Moreover,
  $$E\prq{\psi_n(\oXd^N(s))ds}\leq E\prq{\sigma^N(s)^nds},$$
  hence, due to \eqref{eq:sigma_uniform_bound}, part \ref{eq:moments_bound_sup_bounded} in Assumption \ref{ass:moments_bound} is verified, as well. Assumptions \ref{ass:discrete=fast} to \ref{ass:existence_limit} also hold for the systems with modified rates. Moreover, due to \eqref{eq:m_and_M}, the solution of the deterministic system $(\G_c, \overline{\K}_c)$ coincide with $z$, the solution of the deterministic system $\Sy_c$. Therefore, Assumption \ref{ass:existence_limit} is satisfied as well and we can apply Theorem \ref{thm:weak_convergence} to the modified reaction systems $(\G_c,\overline{\K}^N)$. We have  
  \begin{equation}\label{eq:weak_limit_modified_rates}
   \sup_{t\in[0,T]}\abs{N^{-1}\oXc^N(t)-z(t)}\xrightarrow[N\to\infty]{}0
  \end{equation}
  in probability. Since by definition any function $g\colon\RR^{|\Sp_d|}_{\geq0}\to\RR$ with at most polynomial growth rate satisfies 
  $$\limsup_{\|v\|\to\infty}\frac{|g(v)|}{\psi_n(v)}=0$$
  for $n$ large enough, we have
  $$\sup_{t\in[0,T]}\abs{\int_0^t\Big(g(\oXd^N(s))-E[g(J^{z(s)})]\Big)ds}\xrightarrow[N\to\infty]{}0$$
  in probability. Moreover, for any continuous bounded function $\varphi\colon\RR^{|\Sp_d|}_{\geq0}\to\RR$ and any $\delta\in(0,T)$ we have
  $$\sup_{t\in[\delta,T]}\abs{E[\varphi(\oXd^N(t))]-E[\varphi(J^{z(t)})]}\xrightarrow[N\to\infty]{}0.$$
  The proof is completed by noting that if the path of $\oXd^N$ is different from the path of $\Xd^N$, then we have 
  $$\inf_{S_i\in\Sp_c, t\in[0,T]}N^{-1}\overline{X}^N_i(t)\leq m\quad\text{or}\quad\sup_{S_i\in\Sp_c, t\in[0,T]}N^{-1}\overline{X}_i(t)\geq M.$$
  However, by \eqref{eq:m_and_M} and \eqref{eq:weak_limit_modified_rates}, we have
  \begin{equation*}
  P\pr{\inf_{S_i\in\Sp_c, t\in[0,T]}N^{-1}\overline{X}^N_i(t)\leq m\quad\text{or}\quad\sup_{S_i\in\Sp_c, t\in[0,T]}N^{-1}\overline{X}_i(t)\geq M}\xrightarrow[N\to\infty]{}0.  
  \end{equation*}
  \phantom{.}
 \end{proof}

 \section{ACR setting}
 
 We turn to the ACR setting and start with the formal definition of absolute concentration robustness (ACR).
 
 \begin{definition}\label{def:ACR}
  Let $\Sy=(\G,\K)$ be a reaction system. We say that a species $S_i$ possesses \emph{absolute concentration robustness} (ACR) in $\Sy$ if for any two positive equilibria $q,q'$ of the deterministically modeled system $\Sy$, we have $q_i=q_i'$. In this case, the species $S_i$ is called an \emph{$ACR$ species} and, if a positive equilibrium $q$ exists, $q_i$ is called an \emph{ACR value}. If a system $\Sy$ possesses a non empty set of ACR species, we call it an ACR system. 
 \end{definition}
 
 Consider a reaction system $\Sy$ that has no equilibria or a unique equilibrium. According to Definition \ref{def:ACR}, all the species of $\Sy$ are ACR species, however in these cases the ACR property is not particularly meaningful.
 
 \begin{definition}\label{def:ACR_deg}
  We say that a system $\Sy$ is a \emph{non-degenerate} ACR system if it is an ACR system and possesses at least two positive equilibria. If an ACR system exhibits less than two positive equilibria, we call it a \emph{degenerate} ACR system.
 \end{definition}

We will focus on non-degenerate ACR systems. Note that in such systems not all  species can be ACR species.
 
 In non-degenerate ACR systems, the ACR species maintain their steady-state concentration regardless of the total amount of molecules present in the system. Our goal is to study the behaviour of the system when the abundances of species that do not exhibit ACR tend to infinity. It is therefore natural to use the setting developed in the Section \ref{sec:general_setting} and let the ACR species, or at least a chosen subset of them, be discrete species. We further assume that the rate functions are rescaled consistently with the hypotheses of Section \ref{sec:general_setting}, such that \eqref{eq:lambda_convergence} holds uniformly on compact sets.
 
 In order to study the limiting behaviour of ACR systems, we first introduce the system $\Syl=(\Sp,\C,\fastR,\fastK)$, where the set of reaction $\fastR$ is as defined in \eqref{eq:fastR}, namely
 $$\fastR= \{y\to y'\in\R\,:\,\pi_c(y)\neq0\}.$$
 Furthermore, the kinetics $\fastK$ is given by the functions $\widetilde{\lambda}_r$ defined in \eqref{eq:fast_kinetics}. For convenience, we repeat here the definition: for any $y_r\to y_r'\in\fastR$
  $$\widetilde{\lambda}_{r}(v,w)=\kappa_r(w)v^{\pi_d(y_r)},$$
 where the $\kappa_r(w)$ are as in Assumption \ref{ass:complex_balanced}.  We now state some corollaries of Theorem \ref{thm:weak_convergence}, assuming the next assumption is satisfied.
 
 \begin{assumption}\label{ass:ACR}
  We assume that $\Syl$ is a non-degenerate ACR system, and that at least one of the ACR species is a discrete species.
 \end{assumption}
 
 Consider Example \ref{ex:simple}. The system $\Syl$ coincides with the system introduced in the example itself, namely \eqref{eq:ex_simple}, and it is a non-degenerate ACR system. On the other hand, in Example \ref{ex:non_mass action} the system $\Syl$ is given by
 \begin{gather*}
  A + 2B \xrightarrow{\phantom{\kappa_3}} 3B, \qquad 
  B\ce{<=>[\kappa_1][\kappa_2]} C \xrightarrow{\kappa_3} A,
 \end{gather*}
 with
 $\lambda_{A+2B\to3B}(x)=\kappa_0x_1x_2.$
 Therefore, the dynamics of $\Syl$, if deterministically modeled, coincide with that of the mass action system
 
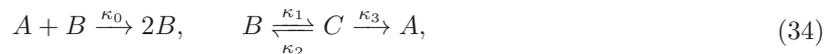
\begin{gather}
  \begin{split}\label{eq:ex_not_mass_action_ACR_system}
   A+B\xrightarrow{\kappa_0} 2B, \qquad
   B\ce{<=>[\kappa_1][\kappa_2]} C \xrightarrow{\kappa_3} A,
  \end{split}
 \end{gather}
 which is a non-degenerate ACR system with equilibria determined by the equations
 $$x_1=\frac{\kappa_1\kappa_3}{\kappa_0(\kappa_2+\kappa_3)}\quad\text{and}\quad x_2=\frac{\kappa_2+\kappa_3}{\kappa_1}x_3.$$
 
 When Assumption \ref{ass:ACR} is fulfilled, denote by $\Spacr$ the set of discrete ACR species. Let $\piacr\colon\RR^{|\Sp_d|}\to\RR^{|\Spacr|}$ be the projection onto the species of $\Spacr$, and let $\Xacr(t)=\piacr(X_d(t))$. Finally, let $q$ be the vector of the ACR values for the species in $\Spacr$, and let 
 \begin{equation}\label{eq:J}
 	J\sim \Pois(q).
 \end{equation}
 \begin{corollary}\label{cor:ACR_equilibrium_point}
 Suppose that Assumptions \ref{ass:discrete=fast}, to \ref{ass:existence_limit} and \ref{ass:ACR} hold. Moreover, assume that \eqref{eq:extra_condition_main_corollary_above} and \eqref{eq:extra_condition_main_corollary_below} hold, and that in the support of any complex $y\in\C$ at most one discrete species appears, and appears with stoichiometric coefficient 1. Let $\pi_c(X_0)$ be a positive equilibrium point for $\Sy_c$, and let $J$ be as in \eqref{eq:J}. Then, for any continuous function $\widehat{g}\colon\RR^{|\Spacr|}_{\geq0}\to\RR$ with at most polynomial growth rate we have
  \begin{equation}\label{eq:statement_cor_ACR_eq_point}
   \sup_{t\in[0,T]}\abs{\int_0^t\Big(\widehat{g}(\Xacr^N(s))-E[\widehat{g}(J)]\Big)ds}\xrightarrow[N\to\infty]{}0
  \end{equation}
  in probability. Moreover, for any continuous bounded function $\widehat{\varphi}\colon\RR^{|\Spacr|}_{\geq0}\to\RR$ and any $\delta\in(0,T)$, we have
  $$\sup_{t\in[\delta,T]}\abs{E[\widehat{\varphi}(\Xacr^N(t))]-E[\widehat{\varphi}(J)]}\xrightarrow[N\to\infty]{}0.$$
 \end{corollary}

The corollary is expressed in standard probabilistic terms, however its meaning might be more intuitively clear with a particular choice of functions $\widehat{g}$ and $\widehat{\varphi}$. For example, if $\widehat{g}$ is taken to be the translated projection $x_i-q_i$, we obtain
$$\sup_{t\in[0,T]}\int_0^t\Big(\Xacr^N(s)-q\Big)ds\xrightarrow[N\to\infty]{}0$$
in probability, where the integral is to be interpreted component-wise. The latter means that on average the counts of the discrete ACR species are well approximated by their ACR value. Moreover, if we let $\widehat{\varphi}$ be  the indicator function $\mathbbm{1}_A$ for a set $A\subseteq\NN^{|\Spacr|}$, we have
$$\sup_{t\in[\delta,T]}\abs{P\pr{\Xacr^N(t)\in A}-P\pr{J\in A}}\xrightarrow[N\to\infty]{}0.$$

 \begin{proof}
  Since $\pi_c(X_0)$ is an equilibrium point for $\Sy_c$, we have $z(t)=\pi_c(X_0)$ for any $t\in[0,T]$. Moreover, by definition $q^w_d$ is the complex balancing equilibrium point of the system $\Sy^w_d$. By the definition of $\fastK$, $(q_d^{\pi_c(X_0)},\pi_c(X_0))$ is a positive equilibrium point for $\Syl$. Hence, $\piacr\pr{q_d^{\pi_c(X_0)}}=q$ and the result follows from Corollary \ref{cor:discrete_species_unary}, applied to the functions $g=\widehat{g}\circ\piacr$ and $\varphi=\widehat{\varphi}\circ\piacr$.
 \end{proof}
 \begin{example}\label{ex:ACR_equilibrium_point}
  Consider the reaction network in \eqref{eq:ex_simple}. $\Sy_c$ is given by \eqref{eq:ex_simple_continuous_system}, for which any non-negative real point is an equilibrium point. We choose a sequence of starting points such that $X_1^N(0)$ is bounded and $N^{-1}X_2^N(0)$ tends to a positive real number. Therefore, the hypotheses of Corollary \ref{cor:ACR_equilibrium_point} are fulfilled for any positive $T$.  In this case there is only one ACR species, namely $A$, and it is the only discrete species. Hence, as $N\to \infty$ the distribution of the counts of species $A$ tends to  a Poisson distribution with mean given by the ACR value $q=\kappa_2/\kappa_1$. The convergence is both on average and at fixed time points, in the sense of Corollary \ref{cor:ACR_equilibrium_point}.
  \hfill $\square$  
 \end{example}
 \begin{corollary}\label{cor:ACR_basin}
 Suppose that Assumptions \ref{ass:discrete=fast}, to \ref{ass:existence_limit} and \ref{ass:ACR} hold. Moreover, assume that \eqref{eq:extra_condition_main_corollary_above} and \eqref{eq:extra_condition_main_corollary_below} hold, and that in the support of any complex $y\in\C$ at most one discrete species appears, and appears with stoichiometric coefficient 1. Assume that $\pi_c(X_0)$ is in the basin of attraction of an equilibrium point of $\Sy_c$. Then, for any continuous function $\widehat{g}\colon\RR^{|\Spacr|}_{\geq0}\to\RR$ with at most polynomial growth rate, $E[\widehat{g}(J^{z(s)})]$ tends to $E[\widehat{g}(J)]$ for $s\to\infty$, where $J$ is as in \eqref{eq:J}. Moreover,
 \begin{equation}
  \sup_{t\in[0,T]}\abs{\int_0^t\Big(\widehat{g}(\Xacr^N(s))-E[\widehat{g}(J^{z(s)})]\Big)ds}\xrightarrow[N\to\infty]{}0
 \end{equation}
 in probability,  where $J^w\sim\Pois(q^w_d)$ and $z(t)$ is as in Assumption \ref{ass:existence_limit}. Finally, for any continuous bounded function $\widehat{\varphi}\colon\RR^{|\Spacr|}_{\geq0}\to\RR$ and any $\delta\in(0,T)$, we have
  $$\sup_{t\in[\delta,T]}\abs{E[\widehat{\varphi}(\Xacr^N(t))]-E[\widehat{\varphi}(J^{z(t)})]}\xrightarrow[N\to\infty]{}0.$$
 \end{corollary}
 The result implies that $E\prq{J^{z(s)}}=q^{z(s)}_d$ tends to $q$ for $s\to\infty$, where $J$ is as in \eqref{eq:J}. Moreover, 
 \begin{equation}
  \sup_{t\in[0,T]}\left|\int_0^t\Big(\Xacr^N(s)-q^{z(s)}_d\Big)ds\right|\xrightarrow[N\to\infty]{}0
 \end{equation}
 in probability. That is, the discrete ACR species are well approximated on average by their ACR value, after a certain time. Finally, for a set $A\subseteq\NN^{|\Spacr|}$
 $$\sup_{t\in[\delta,T]}\abs{P\pr{\Xacr^N(t)\in A}-P\pr{J^{z(t)}\in A}}\xrightarrow[N\to\infty]{}0$$
 and the parameter of $J^{z(t)}$, which is $q^{z(t)}_d$, tends to $q$.
 \begin{proof}
  For the sake of simplicity, throughout this proof $w$ will denote a vector varying in $\RR^{|\Sp_c|}_{>0}$, even if not explicitly stated.
  
  First, recall that, regardless of the value $w>0$, the complex balanced equilibrium $q^w_d$ is the unique solution of a system of multivariate polynomial equations, and as such, it is continuous in the coefficients of the polynomials \cite[Chapter~8]{stetter:numerical}. In particular, it is a continuous function of $w>0$.
  
  Consider a sequence of vectors $(w_n)_{n\in\NN}\subset\RR^{|\Sp_c|}_{>0}$ converging to $w^*>0$. Therefore, the sequence $(q^{w_n}_d)_{n\in\NN}$ converges to $q^{w^*}_d$ and by Lebesgue's Dominated Convergence Theorem we have that
  $$E[\widehat{g}(\piacr(J^{w_n}))]\xrightarrow[n\to\infty]{}E[\widehat{g}(\piacr(J^{w^*}))],$$
  where $J^w\sim\Pois(q^w_d)$. This implies that $E[\widehat{g}(\piacr(J^w))]$ is a continuous function of $w$.
  
  Let $w^*$ be the equilibrium point of $\Sy_c$ whose basin of attraction contains $\pi_c(X_0)$. Since $(q^{w^*}_d,w^*)$ is an equilibrium point of the system $\Syl$ considered in Assumption \ref{ass:ACR}, we have that $\piacr(q^{w^*}_d)=q$. In particular, $\piacr(J^{w^*})\sim J$ and $E[\widehat{g}(J^{z(s)})]$ tends to $E[\widehat{g}(J)]$ for $s\to\infty$. 

  The conclusion of the proof follows directly from Corollary \ref{cor:discrete_species_unary}.
 \end{proof}
 \begin{example}\label{ex:ACR_basin}
  Consider Example \ref{ex:non_mass action}. Let $X_1^N(0)$ be a bounded sequence, and let $N^{-1}(X_2^N(0),X_3^N(0))$ tend to a positive real vector. The continuous system $\Sy_c$, when deterministically modeled, is equivalent to \eqref{eq:ex_not_mass_action_continuous_system}. Therefore the hypotheses of Corollary \ref{cor:ACR_basin} are fulfilled, since any $(b,c)\in\RR_{>0}^2$ is in the basin of attraction of
  $$w^*=\pr{\frac{(\kappa_2+\kappa_3)(b+c)}{\kappa_1+\kappa_2+\kappa_3},\frac{\kappa_1(b+c)}{\kappa_1+\kappa_2+\kappa_3}}.$$
  Let $q$ be the ACR value for $A$ in the system \eqref{eq:ex_not_mass_action_ACR_system}. Therefore, after some time the counts of $A$ are approximately distributed, both on average and at any fixed time point, as a Poisson random variable with mean $q$, in the sense of Corollary \ref{cor:ACR_basin}.
  \hfill $\square$
 \end{example}
 We conclude this section with two examples that do not fit in our theory. However, the examples are still tractable with the techniques we presented.
 \begin{example}\label{ex:not_poisson}
  Consider the stochastic mass action system
  \begin{align*}
   2A+B&\xrightarrow{\kappa_1}3B, \qquad
      B\xrightarrow{\kappa_2}A
  \end{align*}
  Due to \cite{shinar:structural}, we know that $A$ possesses ACR and the system is a non-degenerate ACR system. The discrete system $\Sy^w_d$ is given by
  $$2A\xrightarrow{\kappa_1w}0\xrightarrow{\kappa_2w}A.$$
  Therefore, the discrete system is not weakly reversible and by Theorem \ref{thm:complex_balanced} it cannot be complex balanced. It follows that Assumption \ref{ass:complex_balanced} does not hold, and by \cite{cappelletti:complex_balanced} we know that the $\Sy^w_d$, stochastically modeled, cannot exhibit a Poisson stationary distribution. However, a stationary distribution $\mu^w$ can be shown to exist by standard methods. By following the proofs of the results in this paper, we can still argue that, if $X^N(0)=(a,bN)$ and $J^w\sim\mu^w$, for any continuous function $g\colon\RR_{\geq0}\to\RR$ with at most polynomial growth rate we have 
  \begin{equation}
   \sup_{t\in[0,T]}\abs{\int_0^t\Big(g(X_1^N(s))-E[g(J^{z(s)})]\Big)ds}\xrightarrow[N\to\infty]{}0
  \end{equation}
  in probability. Moreover, for any continuous bounded function $\varphi\colon\RR_{\geq0}\to\RR$ and any $\delta\in(0,T)$, we have
  $$\sup_{t\in[\delta,T]}\abs{E[\varphi(X_1^N(t))]-E[\varphi(J^{z(t)})]}\xrightarrow[N\to\infty]{}0.$$
  In this regard, see also Remark \ref{rem:general_stationary_distribution}, where this situation is discussed.

  Unfortunately, there are not many methods available that explicitly calculate stationary distributions for non-complex balanced systems.  Thus, most examples not fulfilling Assumption \ref{ass:complex_balanced} are analytically intractable. However, for some calculations of  stationary distributions for non-complex balanced systems, see \cite{anderson:lyapunov}.
  \hfill $\square$
 \end{example}
 
 \begin{example}\label{ex:bimolecular}
  The results in this section rely on Corollary \ref{cor:discrete_species_unary}, whose assumptions are easy to check. However, here we show how the same conclusions can be derived directly from Theorem \ref{thm:weak_convergence}, provided that Assumption \ref{ass:moments_bound} can be checked for the process or for a convenient modification thereof. Consider the deterministic mass action system given by
  
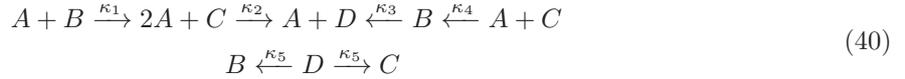
\begin{align}
  \begin{split}
   A+B \ce{->[\kappa_1]} 2A+C \ce{->[\kappa_2]} A&+D \ce{<-[\kappa_3]} B \ce{<-[\kappa_4]} A+C\\
   B\ce{<-[\kappa_5]} &D \ce{->[\kappa_5]} C
   \end{split}
   \label{eq:856786567}
  \end{align}
  with the constraint that $\kappa_1/\kappa_2=\kappa_3/\kappa_4$. The species $A$ is the only ACR species, and its ACR value is
  $$q=\sqrt{\frac{\kappa_3\kappa_4}{\kappa_1\kappa_2}}=\frac{\kappa_4}{\kappa_2},$$
   where the last equality derives from $\kappa_1/\kappa_2=\kappa_3/\kappa_4$.  If we let $A$ be the discrete species and $B$, $C$ and $D$ be the continuous species, then Assumption \ref{ass:discrete=fast} is fulfilled. In this case, the system $\Syl$ coincides with the above mass action system \eqref{eq:856786567}, and Assumption \ref{ass:ACR} is also satisfied. Let $w=(w_1,w_2,w_3)$ denote the concentrations of the species $B$, $C$ and $D$ respectively. Then, the system $\Sy^w_d$ is given by the stochastic mass action system
  $$0\ce{<=>[\kappa_3w_1][\kappa_4w_2]}A\ce{<=>[\kappa_1w_1][\kappa_2w_2]}2A,$$
  which is complex balanced for any positive $w$ due to the assumption $\kappa_1/\kappa_2=\kappa_3/\kappa_4$. The system $\Sy^w_d$ is also irreducible. Therefore, Assumption \ref{ass:complex_balanced} is fulfilled and the complex balanced equilibrium is given by
  $$q^w_d=\frac{\kappa_1w_1}{\kappa_2w_2}.$$
  The system $\Sy_c$ is given by
  $$
   B \ce{<=>} C \ce{<=>[][\kappa_5]}D\ce{<=>[\kappa_5][\kappa_3]}B,
  $$
  where
  \begin{align*}
   \lambda_{B\to C}(w)= \lambda_{C\to D}(w)= \frac{\kappa_1^2w_1^2}{\kappa_2w_2}\quad \text{and} \quad 
   \lambda_{C\to B}(w)=\frac{\kappa_1\kappa_4w_1}{\kappa_2}.
  \end{align*}
  For simplicity, we choose our initial condition equal to an equilibrium of \eqref{eq:856786567}
 $$X^N(0)=\pr{\frac{\kappa_4}{\kappa_2}, N, \frac{\kappa_1}{\kappa_4}N, \frac{\kappa_3}{\kappa_5}N}.$$
 In this case, we have that $z(t)$ is constantly equal to $z(0)=(1,\kappa_1/\kappa_4, \kappa_3/\kappa_5)$ and Assumption \ref{ass:existence_limit} is clearly satisfied.
 
 Following the proof of Corollary \ref{cor:discrete_species_unary}, we can modify the kinetics of the original system by
$$\overline{\lambda}^N_r(x)=\lambda^N_r\pr{\pi^N(x)}\quad\text{for }x\in\NN^{|\Sp|}.$$
We can derive the same conclusions as in Corollary \ref{cor:discrete_species_unary} by noting that $\Sy^w_d$ defines a birth-death process whose moments are all uniformly bounded in time, and a uniform bound can be found also by letting $w$ vary in a compact set. Therefore,  the counts of the ACR species $A$ are well approximated by a Poisson random variable with mean the ACR value $\kappa_4/\kappa_2$. The approximation is both on average and at finite time intervals, in the sense of Corollary \ref{cor:ACR_equilibrium_point}. If, on the other hand, $\pi_c(X^N(0))$ does not converge to an equilibrium point of $\Sy_c$, then we can still control the behaviour of the species $A$ in the sense of Corollary \ref{cor:ACR_basin}.
 
 \end{example}
 
 \section{EnvZ/OmpR signaling system}
 
 As another application of our results, we consider the two-component EnvZ/OmpR osmoregulatory signaling system in \textit{Escherichia coli}, using the model proposed in \cite{shinar:input} and considered in \cite{anderson:ACR,shinar:structural}. The model corresponds to the following mass action system, where $[D]$ and $[T]$ are constants describing the concentration of some slowly interacting chemical species:
 \begin{gather*}
  XD\ce{<=>[\kappa_1][\kappa_2{[D]}]}X\ce{<=>[\kappa_3{[T]}][\kappa_4]}XT\ce{->[\kappa_5]}X_p\\
  X_p+Y\ce{<=>[\kappa_6][\kappa_7]}X_pY\ce{->[\kappa_8]}X+Y_p\\
  XD+Y_p\ce{<=>[\kappa_9][\kappa_{10}]}XDY_p\ce{->[\kappa_{11}]}XD+Y\;,
 \end{gather*}
 where $X=\text{EnvZ}, Y=\text{OmpR}, X_p= \text{EnvZ-P}, Y_p=\text{OmpR-P}, D=\text{ADP}$, and $T=\text{ATP}$. The abundances of both ADP and ATP are assumed to be  large enough so that their consumption in the first chain of reactions only negligibly changes their concentration. The first chain of reactions describes the phosphorylization of EnvZ, the second chain corresponds to the transfer of the phosphate group from EnvZ to OmpR, and finally the third chain describes the dephosphorylization of OmpR.

 Due to \cite{anderson:ACR,shinar:structural}, it is known that the species $Y_p$ exhibits ACR. Moreover, $Y_p$ is the only ACR species. For simplicity, instead of ordering the species, here we will denote by $z_S$ the concentration of the species $S$. At equilibrium
 \begin{align*}
  z_{Y_p}&=\frac{\kappa_1\kappa_3\kappa_5(\kappa_{10}+\kappa_{11})[T]}{\kappa_2(\kappa_4+\kappa_5)\kappa_9\kappa_{11}[D]}=q\\
  z_{XT}&=\frac{\kappa_3[T]}{\kappa_4+\kappa_5}z_X=\frac{\kappa_1}{\kappa_4+\kappa_5}z_{XD}=\frac{\kappa_8}{\kappa_5}z_{X_pY}=\frac{\kappa_{11}}{\kappa_5}z_{XDY_p}\\
  z_{X_p}z_Y&=\frac{\kappa_7}{\kappa_6}z_{X_pY}+\frac{\kappa_{11}}{\kappa_6}z_{XDY_p}.
 \end{align*}
 There are two quantities that are conserved at any time point, namely
 \begin{align*}
  c_1&=z_Y(t)+z_{Y_p}(t)+z_{X_pY}(t)+z_{XDY_p}(t)\\
  c_2&=z_{X_p}(t)+z_{XT}(t)+z_X(t)+z_{XD}(t)+z_{X_pY}(t)+z_{XDY_p}(t),
 \end{align*}
 for some positive constants $c_1,c_2$ depending on the initial conditions. If the amounts $c_1$ and $c_2$ are increased, then the equilibrium concentrations of all the species not exhibiting ACR are increased as well, except for $X_p$ and $Y$, the equilibrium concentration of one of which could remain small.
 
 Consider now the above reaction system in the stochastic setting. We want to know what happens if we increase the initial counts of the species such that the conserved amounts are equally increased and the initial condition is in a neighbourhood of an equilibrium point of the system. Therefore, we uniformly increase the counts of the species not exhibiting ACR, and we choose to keep $X_p$ or $Y$ small. We consider a sequence of processes $X^N$ indexed by $N\in\NN$, which are associated with the above reaction system. We assume that $X^N(0)$ is such that the entries relative to $Y_p$ and $Y$, denoted by $X^N_{Y_p}(0)$ and $X^N_{Y}(0)$ respectively, are bounded by a constant $B$, and that all the other entries, if rescaled by $N$, converge to some positive number.
 In this setting, the discrete species are $Y$ and $Y_p$, we have $\beta_r=1$ for any reaction $y_r\to y_r'$ of the system and Assumption \ref{ass:discrete=fast} is fulfilled. For any positive vector $w$ of continuous species concentrations, the system $\Sy^w_d$ is given by
 \begin{equation}\label{eq:network}
0 \ce{<=>[{\kappa_7 w_{X_pY}+\kappa_{11} w_{XDY_p}}][{\kappa_6 w_{X_p}}]}Y\qquad 0\ce{<=>[{\kappa_8 w_{X_pY} + \kappa_{10} w_{XDY_p}}][{\kappa_9 w_{XD}}]}Y_p,  
 \end{equation}
 and Assumption \ref{ass:complex_balanced} holds thanks to Theorem \ref{thm:deficiency_zero_iff}. It is not difficult to check that the mass action system \eqref{eq:network} is irreducible, since every state of possible counts of $Y$ and $Y_p$ is accessible from any other state. The complex balanced equilibrium $q^w_d$ is given by
 $$q^w_d=\pr{\frac{\kappa_7 w_{X_pY}+\kappa_{11} w_{XDY_p}}{\kappa_6 w_{X_p}},\frac{\kappa_8 w_{X_pY} + \kappa_{10} w_{XDY_p}}{\kappa_9 w_{XD}}},$$
 where the first entry refers to $Y$ and the second one to $Y_p$. The system $\Sy_c$ is given by
 \begin{gather*}
  XD\ce{<=>[\kappa_1][\kappa_2{[D]}]}X\ce{<=>[\kappa_3{[T]}][\kappa_4]}XT\ce{->[\kappa_5]}X_p\\
  X_p\ce{<=>[][\kappa_7]}X_pY\ce{->[\kappa_8]}X\\
  XD\ce{<=>[][\kappa_{10}]}XDY_p\ce{->[\kappa_{11}]}XD\;,
 \end{gather*}
 with
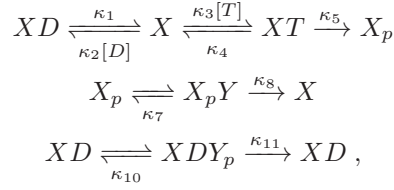
 \begin{align*}
  \lambda_{X_p\to X_pY}(w)&=\kappa_7 w_{X_pY}+\kappa_{11} w_{XDY_p}\\
  \lambda_{XD\to XDY_p}(w)&=\kappa_8 w_{X_pY} + \kappa_{10} w_{XDY_p}.
 \end{align*}
 The equilibria of the system are the positive vectors $w^*$ that satisfy
 $$w^*_{XT}=\frac{\kappa_3[T]}{\kappa_4+\kappa_5}w^*_X=\frac{\kappa_1}{\kappa_4+\kappa_5}w^*_{XD}=\frac{\kappa_8}{\kappa_5}w^*_{X_pY}=\frac{\kappa_{11}}{\kappa_5}w^*_{XDY_p}.$$
 If $\pi_c(X^N(0))$ is such a vector $w^*$, or belongs to its basin of attraction, then Assumption \ref{ass:existence_limit} holds. Assumption \ref{ass:ACR} also holds, since the system $\Syl$ corresponds to the original EnvZ/OmpR signaling system, and the unique ACR species $Y_p$ is discrete. By making use of the fact that the original system is mass action kinetics and by Remark \ref{rem:mass_action}, it is easy to see that the remaining assumptions of Corollary \ref{cor:ACR_equilibrium_point} (if $\pi_c(X^N(0))$ is an equilibrium $w^*$) or of Corollary \ref{cor:ACR_basin} (if $\pi_c(X^N(0))$ is in the basin of attraction of an equilibrium $w^*$) are fulfilled, and the results can be applied. Therefore, $X^N_{Y_p}(t)$ can be approximated by a Poisson random variable $J$ with mean $q$, both on average and at any fixed time point, in the sense of Corollary \ref{cor:ACR_equilibrium_point} or Corollary \ref{cor:ACR_basin}. The results are in accordance with the simulations in \cite{anderson:ACR}. 
 
 Alternatively, we could have applied the results of this paper to the signaling network by considering $Y_p$ as the only discrete species, therefore increasing the initial counts of all other species, and by letting $X^N$ be the process associated with
 \begin{gather*}
  XD\ce{<=>[\kappa_1][\kappa_2{[D]}]}X\ce{<=>[\kappa_3{[T]}][\kappa_4]}XT\ce{->[\kappa_5]}X_p\\
  X_p+Y\ce{<=>[\kappa_6/N][\kappa_7]}X_pY\ce{->[\kappa_8]}X+Y_p\\
  XD+Y_p\ce{<=>[\kappa_9][\kappa_{10}]}XDY_p\ce{->[\kappa_{11}]}XD+Y\;,
 \end{gather*}
 where $\kappa_6$ has been rescaled. With this choice of rescaling, due to Remark \ref{rem:rate_constants_convergence}, we have that $\beta_r=1$ for any reaction $y_r\to y_r'$ of the system. Our results can be used to draw the same conclusion as before in this different setting.
 
  \section*{Acknowledgements}
 
 The authors would like to thank Carsten Wiuf: this paper is the result of a collaboration among the authors which has been encouraged and made possible by him, as part of Cappelletti's Ph.D.\ program.
 
 \bibliographystyle{spmpsci}
 \bibliography{bib}
\end{document}